\documentclass[10pt, reqno]{amsart}
\usepackage{amsmath,amssymb,amsbsy,amsfonts,latexsym,amsthm,amscd,amsxtra, xfrac,
amssymb}
\usepackage[all]{xy}

\newtheorem{Teo}{Theorem}
\newtheorem{Def}{Definition}
\newtheorem{Prop}{Proposition}
\newtheorem{Cor}{Corollary}

\newtheorem{Rema}{Remark}
\newenvironment{Rem}{\begin{Rema} \begin{upshape}} {\end{upshape}\end{Rema}}
\newtheorem*{Pf}{Proof}
\newenvironment{Proof}{\begin{Pf} \begin{upshape}} {\end{upshape} \qed\end{Pf}}

\newcommand\beqa[1]{ \begin{eqnarray} \label{#1}}
\newcommand{\eeqa}{ \end{eqnarray} }
\newcommand{\beqano}{ \begin{eqnarray*} }
\newcommand{\eeqano}{ \end{eqnarray*} }

\newcommand{\T}{ {\mathbb T}   }

\newcommand{\R}{ {\mathbb R}   }
\newcommand{\Z}{ {\mathbb Z}   }

\newcommand{\Q}{ {\mathbb Q}   }

\newcommand{\cD}{ {\mathcal D}   }

\def\beal{\begin{aligned}}
\def\enal{\end{aligned}}

\renewcommand \a {\alpha}
\newcommand \e {\varepsilon }

\renewcommand \b  {\beta}

\newcommand \D{\Delta}
\newcommand \m {\mu}

\newcommand \Om {\Omega}

\newcommand \f {\varphi}

\newcommand \g {\gamma}

\newcommand \G {\Gamma}

\renewcommand \th {\theta}

\renewcommand \phi {\varphi}

\newcommand \cE {{\mathcal E}}

\def\ie{\hbox{\it i.e.\ }}

\def\~{\tilde}
\def\Bbb{\mathbb}

\def\A{\Bbb A}

\def\R{\Bbb R}

\def\Z{\Bbb Z}
\def\T{\Bbb T}

\def\Q{\mathbb Q}

\def\I{\mathcal I}
\def\dt{\delta}

\def\th{\theta}

\def\I{\mathcal I}

\begin{document}

\title{Computing Mather's $\beta$-function for Birkhoff billiards}
\author{Alfonso Sorrentino}
\email{sorrentino@mat.uniroma3.it}
\address{Dipartimento di Matematica e Fisica, Sezione di Matematica, Universit\`a degli Studi Roma Tre, Largo S. Leonardo Murialdo 1, 00146 Rome (Italy).}
\date{\today}

\maketitle

\begin{abstract}
This article is concerned with the study of {\em Mather's $\beta$-function}  associated to Birkhoff billiards.
This function corresponds to the minimal average action of orbits with a prescribed rotation number and, from a different perspective, it can be related to the maximal perimeter of periodic orbits with a given rotation number, the so-called {\em Marked length spetrum}.
After having recalled its main properties and its relevance to the study of the billiard dynamics, we  stress its connections to some intriguing open questions:  Birkhoff conjecture and 
the isospectral rigidity of convex billiards. Both these problems, in fact, can be conveniently translated into questions on this function. This motivates our investigation aiming at understanding its main features and properties.
In particular, we provide an explicit representation of the coefficients of its (formal) Taylor expansion at zero, only in terms  of the curvature of the boundary. In the case of integrable billiards, this result provides a representation formula for the $\beta$-function near $0$. Moreover, we apply and check these results in the case of circular and elliptic billiards.
\end{abstract}


\section{Introduction }

In this note we would like to provide explicit computations for {\it Mather's $\beta$-function} (or {\it minimal average action}) in the case of Birkhoff billiards. In particular, we aim at describing an explicit representation of the coefficients of its (formal) Taylor expansion, in terms of the curvature of the boundary. This function 
-- which is related, at least in the case of rational rotation numbers, to the maximal length of periodic orbits with a given rotation number (the so-called {\it marked lenght spetrum}) -- plays a crucial r\^ole
in the comprehension of different rigidity phenomena that appear in the study of convex billiards; moreover, many intriguing unanswered questions and conjectures can be easily translated into questions on this function. Hence, we believe that understanding its main features and properties -- besides being interesting {\it per se} -- is an essention step in order to  tackle and unravel these compelling open questions.\\

A {\it Birkhoff billiard} \footnote{
This conceptually simple model, yet dynamically very rich, has been first introduced
by G. D. Birkhoff \cite{Birkhoff} as a mathematical playground to prove, with as little
technicality as possible, some dynamical applications of Poincare's last geometric
theorem and its generalisations:
\begin{quote}
``[...]{\it This example is very illuminating for the following reason: Any dynamical
system with two degrees of freedom is isomorphic with the motion of a particle on
a smooth surface rotating uniformly about a fixed axis and carrying a conservative
field of force with it} (see \cite{Birkhoff1917}). {\it In particular if the surface
is not rotating and if the field of force is lacking, the paths of the particles
will be geodesics. If the surface is conceived of as convex to begin with and then
gradually to be flattened to the form of a plane convex curve $C$, the `billiard ball'
problems results. But in this problem the formal side, usually so formidable in dynamics,
almost completely disappears, and only the interesting qualitative questions need
to be considered.}[...] ''\\
(G. D. Birkhoff, \cite[pp. 155-156]{Birkhoff})\\
\end{quote}
}
is a dynamical model describing the motion of a mass point
inside a (strictly) convex domain $\Omega \subset \R^2$ with smooth boundary.
The massless {billiard ball} moves with unit velocity and without friction
following a rectilinear path;  when it hits the boundary it reflects {elastically}
according to the standard {\it reflection law}: the angle of reflection is equal to
the angle of incidence. Such trajectories are sometimes called {\it broken geodesics}.

Let us recall some properties of the billiard map. We refer to \cite{Siburg, Tabach} for a more comprehensive introduction to the study of billiards.

\begin{figure} [h!]
\begin{center}
\includegraphics[scale=0.35]{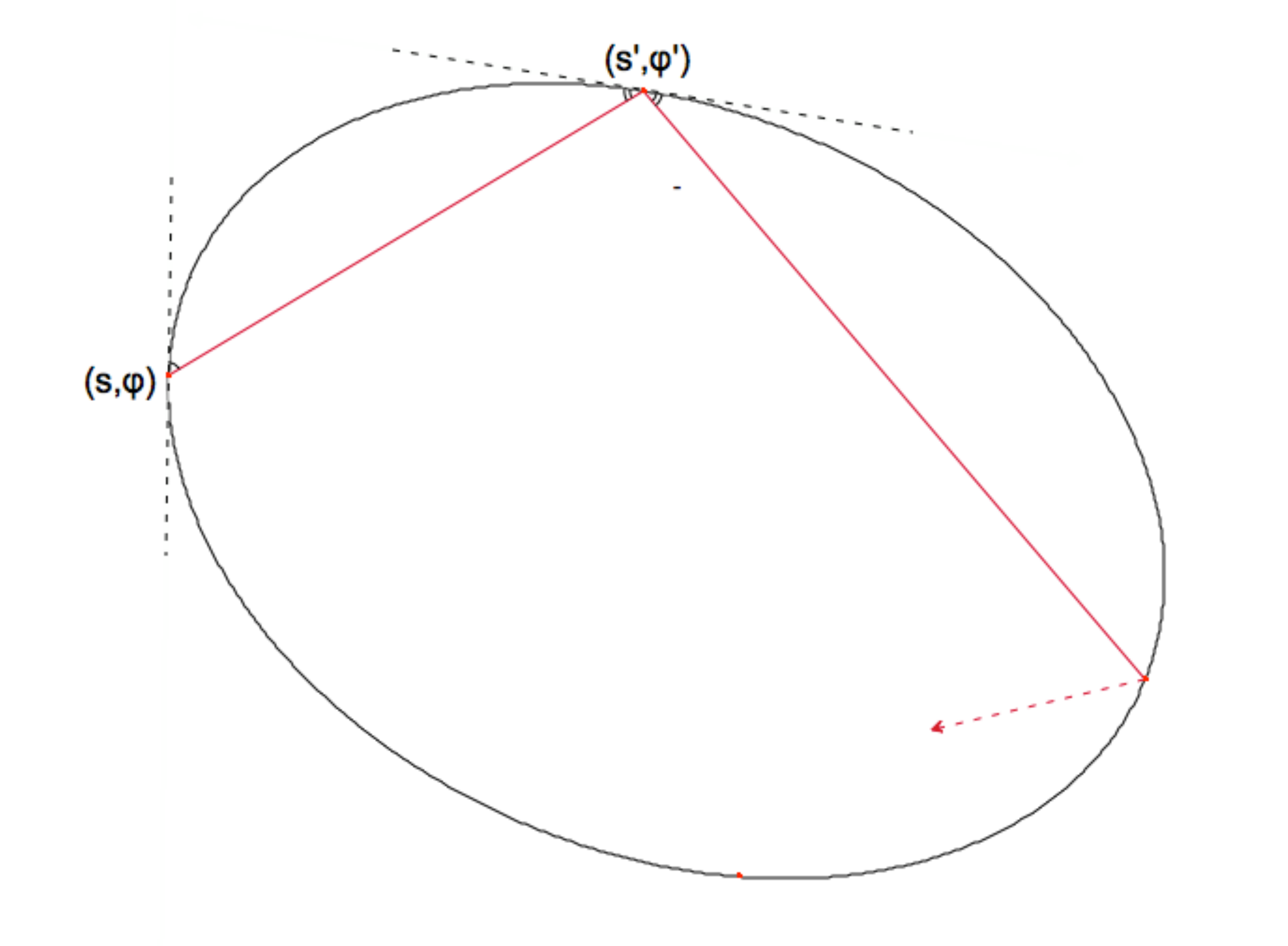}
\caption{}
\label{billiard}
\end{center}
\end{figure}

Let $\Omega$ be a strictly convex domain in $\R^2$ with $C^r$ boundary $\partial \Omega$,
with $r\geq 3$. The phase space $M$ of the billiard map consists of unit vectors
$(x,v)$ whose foot points $x$ are on $\partial \Omega$ and which have inward directions.
The billiard ball map $f:M \longrightarrow M$ takes $(x,v)$ to $(x',v')$, where $x'$
represents the point where the trajectory starting at $x$ with velocity $v$ hits the boundary
$\partial \Omega$ again, and $v'$ is the {\it reflected velocity}, according to
the standard reflection law: angle of incidence is equal to the angle of reflection (figure \ref{billiard}).

\begin{Rem}
Observe that if $\Omega$ is not convex, then the billiard map is not continuous.
Moreover, as pointed out by Halpern \cite{Halpern}, if the boundary is not at
least $C^3$, then the flow might not be complete.
\end{Rem}

Let us introduce coordinates on $M$.
We suppose that $\partial \Omega$ is parametrized  by  arc-length $s$ and
let $\g:  [0, l] \longrightarrow \R^2$ denote such a parametrization,
where $l=l(\partial \Omega)$ denotes the length of $\partial \Omega$. Let $\phi$
be the angle between $v$ and the positive tangent to $\partial \Omega$ at $x$.
Hence, $ M$ can be identified with the annulus $\A = [0,l] \times (0,\pi)$
and the billiard map $f$ can be described as

\begin{eqnarray*}
f: [0,l] \times (0,\pi) &\longrightarrow& [0,l] \times (0,\pi)\\
(s,\phi) &\longmapsto & (s',\phi').
\end{eqnarray*}

In particular $f$ can be extended to $\bar{\A}=[0,l] \times [0,\pi]$ by fixing
$f(s,0)=f(s,\pi)= {\rm Id}$,
for all $s$. \\

Let us denote by
$$
\ell(s,s') := \|\g(s) - \g(s')\|
$$
the Euclidean distance between two points on $\partial \Omega$. It is easy to prove that
\begin{equation}\label{genfunctbill}
\left\{ \begin{array}{l}
\dfrac{\partial \ell}{\partial s}(s,s') = - \cos \phi \\
\\
\dfrac{\partial \ell}{\partial s'}(s,s') = \cos \phi'\,.\\
\end{array}\right.
\end{equation}

\begin{Rem}
If we lift everything to the universal cover and introduce new coordinates
$(x,y)=(s, -\cos \phi) \in \R \times (-1,1)$, then the billiard map is a twist map
with $\ell$ as generating function and it preserves the area form $dx \wedge dy$. See \cite{Siburg, Tabach}.\\
\end{Rem}

Despite the apparently simple (local) dynamics, the qualitative dynamical properties of billiard maps are extremely {non-local}. This global influence
on the dynamics translates into several intriguing {\it rigidity  phenomena}, which
are at the basis of several unanswered questions and conjectures. Amongst many, two noteworthy ones regard the rigidity of the {\it length spectrum} (see subsection 1.1) and the classification of {\it integrable billiards}, also known as {\it Birkhoff conjecture} (see subsection 1.2).
Both questions are deeply tangled to properties of  Mather's $\beta$-function (see definition \ref{defbeta}) and can be translated into questions on its rigidity and regularity, as we shall explain in the following (see subsection 1.3).\\

\noindent  {\bf 1.1 -  Periodic orbits and Marked length spectrum}.

 The study of periodic orbits and their properties have been amongst the first dynamical features of billiards that have been investigated.
One of the first results in the theory of billiards, for example, can be considered  Birkhoff's application of Poincare's last geometric theorem to show
the existence of infinitely many {\it distinct} periodic orbits  \cite{Birkhoff}. Since then, new phenomena have been pointed out and many interesting questions have been raised.

How do we distinguish {\it distinct} periodic orbits? One could try to classify them in terms of their {\it period},  {\it i.e.},  the minimal number of  times that the ball reflects before going back to the initial position with the initial direction. However, while in some cases this quantity allows one to distinguish different periodic orbits, in many cases  it is not sufficient anymore: periodic orbits with the same periods may wind a different  number of times  before closing; this will clearly translate into a different topological shape.

A better invariant that one should consider  is the so-called
{\it rotation number}. The rotation number  of a periodic billiard trajectory
(respectively, a closed broken geodesic) is a rational
number
\[
\dfrac{p}{q}
\ =\
\dfrac{\text{winding number}}
{\text{number of reflections}}\ \in\ \big(0,\frac 12\Big],
\]
where the winding number $p>1$ is defined as follows.
Fix the positive orientation of $\partial \Om$ and
pick any reflection point of the closed geodesic on 
$\partial \Om$; then follow the trajectory and count
how many times it goes around $\partial \Om$ in
the positive direction until it comes back to
the starting point.
Notice that inverting the direction of motion for every
periodic billiard trajectory  of rotation number
$p/q  \in (0, 1/2]$, we obtain a trajectory with rotation number
$(q-p)/q \in [1/2,1)$. \\

In \cite{Birkhoff}, Birkhoff proved that for every $p/q \in (0, 1/2]$ in lowest terms,
there are at least two closed orbits  of rotation number
$p/q$: one maximizing the total length and the other obtained by min-max methods (see also  \cite[Theorem 1.2.4]{Siburg}).
This result is clearly optimal: in the case of a billiard in an ellipse, for example, there are only two periodic orbits of period $2$ (also called {\it diameters}), which correspond to the two semi-axis of the ellipse (see for example subsection 1.2 or Section 3.2). However, it is easy to find cases in which there are more than two periodic orbits for any given rotation number: think, for example, of a billiard in a disk where, due to the existence of a $1$-dimensional group of symmetries (rotations), each periodic orbit generates a $1$-dimensional family of similar ones; for example, all diameters are periodic orbits with period $2$ (see subsection 1.2 and Section 3.1).\\

This raises this natural question:\\

\noindent  {\it What information on the geometry of the billiard domain do closed orbits carry? Does the knowledge of the lengths of periodic orbits allow one to reconstruct the billiard domain? }\\


One could `organize' this set of information in a more functional way, for instance by associating to each length the corresponding rotation number or even refining it by considering only orbits with maximal length amongst those with a given rotation number.
This map is called the ({\it maximal}) {\it marked length spectrum} of $\Omega$.\\

\begin{Def}[{\bf Marked Length Spectrum}] {\it
Given $\Omega$ a strictly convex planar domain with smooth boundary, we define its Marked length spectrum
${\mathcal ML}_\Om: \Q\cap\big(0,\frac 12\Big] \longrightarrow \R_+ $ as:
\begin{eqnarray*}{\mathcal ML}_\Om({p}/{q}) =   \max \Big \{ \mbox{lengths of periodic orbits with rotation number}\; p/q \Big \}.\\
\end{eqnarray*}
}
\end{Def}

\noindent{\bf Question I  (Guillemin--Melrose \cite{GM}).}  {\it Let $\Omega_1$ and $\Omega_2$ be two strictly convex planar domains with smooth boundaries and assume that they are isospectral, {\it i.e.}, 
${\mathcal ML}_{\Om_1} \equiv {\mathcal ML}_{\Om_2}$. Is it true that $\Omega_1$ and $\Omega_2$ are isometric?}\\

\begin{Rem}
The above question could be reformulated -- and it remains still meaningful and interesting -- by asking that they two domains are `only' isospectral near the boundary, {\it i.e.}, ${\mathcal ML}_{\Om_1}(p/q) = {\mathcal ML}_{\Om_2} (p/q)$ for all
$p/q \in \Q \cap [0,\e)$, for some $0<\e\leq 1/2$.\\
\end{Rem}

See subsection 1.3 for a reformulation of this question in terms of Mather's $\beta$ function (Questions I bis and ter).\\



\noindent {\bf 1.2 - Integrable billiards and Birkhoff conjecture}.

\begin{figure} [h!]
\begin{center}
\includegraphics[scale=0.3]{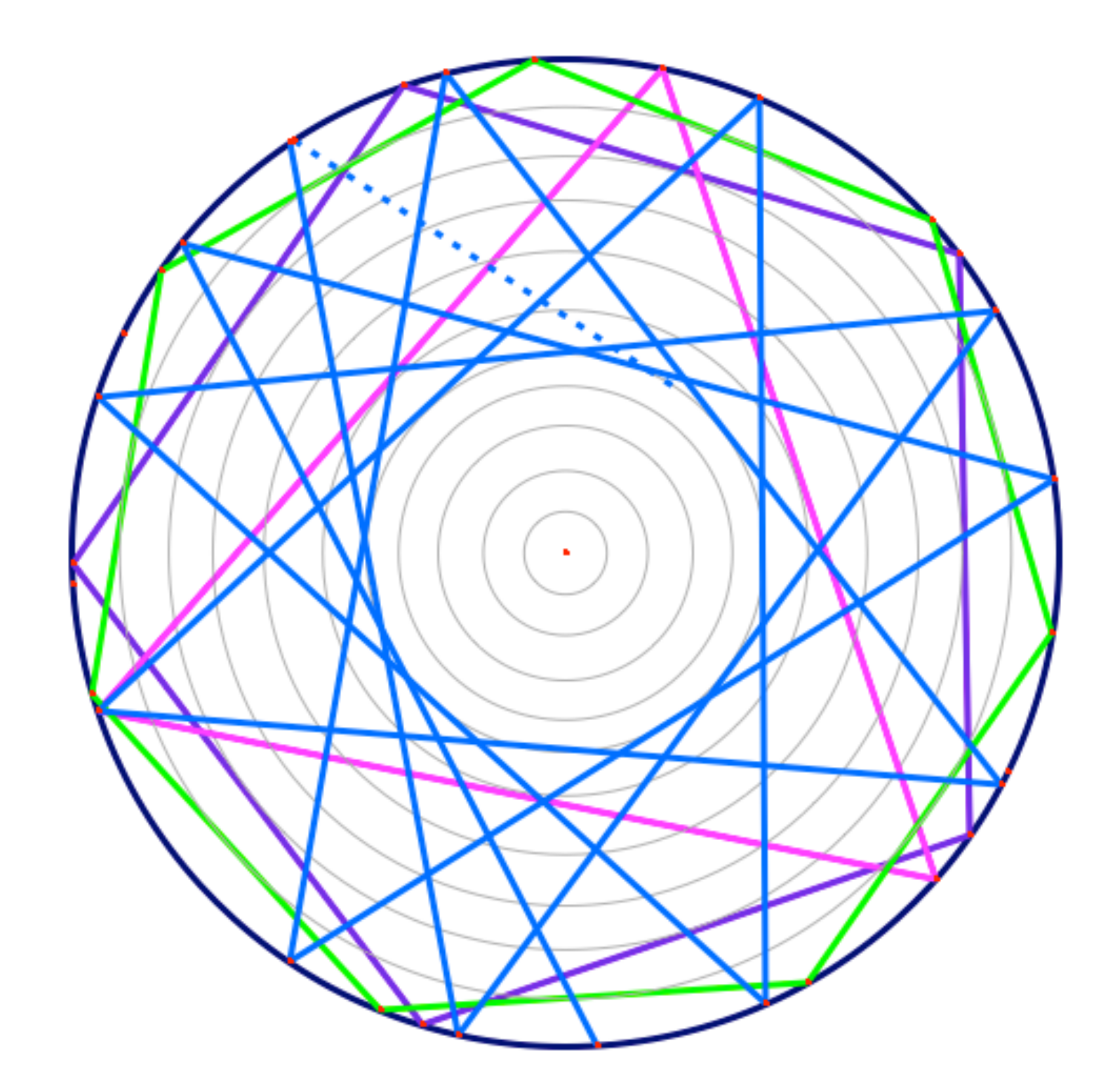}
\caption{Billiard in a disc}
\label{circle-billiard}
\end{center}
\end{figure}

The easiest example of billiard is given by a billiard in a disc $\cD$ (for example of radius $R$). It is easy to check in this case that the angle of reflection remains constant at each reflection (see also  \cite[Chapter 2]{Tabach} and Section \ref{exampledisc}). 
If we denote by $s$  the arc-length parameter ({\it i.e.}, $s\in {\R}/{ {\small 2\pi R \Z}}$) and by $\theta \in (0,\pi/2]$ the angle of reflection, then
the billiard map  has a very simple form:
$$
f(s,\theta) = (s + 2R\, \theta,\; \theta).
$$
In particular, $\theta$ stays constant along the orbit and it represents an {\it integral of motion} for the map.
Moreover, this billiard enjoys the peculiar property of
having  the phase space -- which is topologically a cylinder --  completely foliated by homotopically non-trivial invariant curves ${\mathcal C}_{\theta_0}=\{\theta\equiv \theta_0\}$. These curves correspond to concentric circles of radii 
$\rho_0= R\cos \theta_0$ and are examples of what are called {\it caustics},  {\it i.e.}, (smooth and convex) curves with the property
that if a trajectory is tangent to one of them, then it will remain tangent after each reflection (see figure \ref{circle-billiard}).

A billiard in a disc is an example of an  {\it integrable billiard}. There are different  ways to define global/local integrability for billiards (the equivalence
of these notions is an interesting problem itself):
\begin{itemize}
\item[-] either through the existence of an integral of motion, globally or locally near the boundary (in the circular case  an integral of motion is given by $I(s,\theta)=\theta$),
\item[-] or through the existence of a (smooth) foliation of the whole phase space (or locally in a neighbourhood of the boundary $\{\theta=0\}$), consisting of invariant curves of the billiard map; for example, in the circular case these are given by ${\mathcal C}_{\theta}$. This property translates (under suitable assumptions) into the existence of a (smooth) family of caustics, globally or locally near the boundary (in the circular case, the concentric circles of radii $R\cos \theta$).\\
\end{itemize}


In \cite{Bialy},  Misha Bialy proved the following beautiful result concerning global integrability (see also \cite{Woi}):\\

\noindent {\bf Theorem (Bialy).}
\noindent {\it If the phase space of the billiard ball map is globally foliated by continuous
invariant curves which are not null-homotopic, then  it is a circular billiard.}\\

\begin{figure} [h!]
\begin{center}
\includegraphics[scale=0.2]{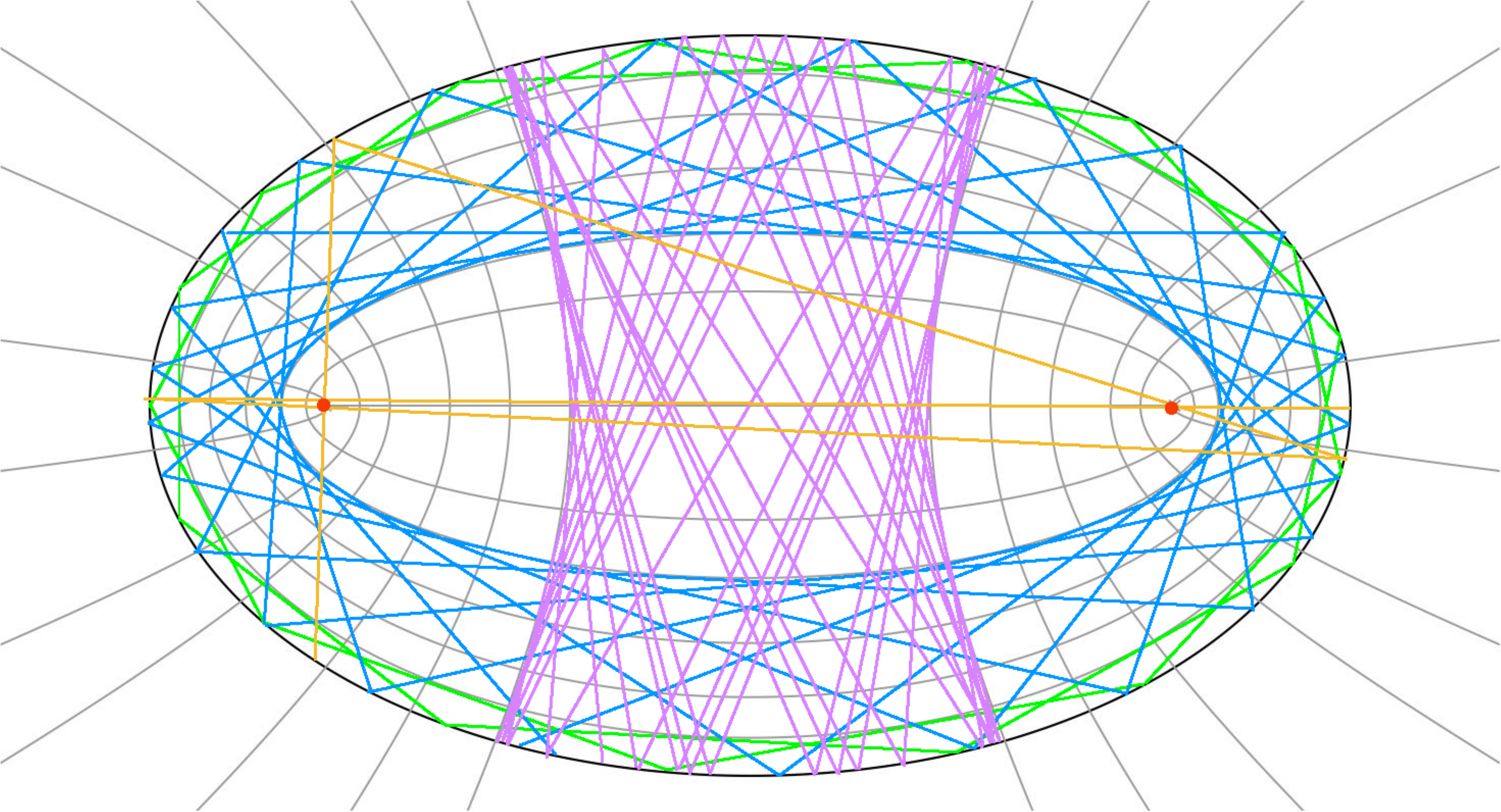}
\caption{Billiard in an ellipse}
\label{ellipse-billiard}
\end{center}
\end{figure}

However, while circular billiards are the only examples of global integrable billiards, local integrability is still an intriguing open question.
One could  consider a  billiard in an ellipse: this is in fact (locally) integrable (see Section 3.2). Yet, the dynamical picture
is very distinct from the circular case: as it is showed in figure \ref{ellipse-billiard}, each trajectory which does not pass through
a focal point, is always tangent to precisely one confocal conic section, either
a confocal ellipse or the two branches of a confocal hyperbola (see for example
\cite[Chapter 4]{Tabach}). Thus, the confocal ellipses inside an elliptical billiards
are convex caustics, but they do not foliate the whole domain: the segment between
the two foci is left out (describing the dynamics explicitly is much more complicated: see for example \cite{Taba} and Section \ref{exampleellipse}). \\

\noindent{\bf Question II (Birkhoff).} {\it Are there other examples of (locally) integrable billiards?}\\

\noindent A negative answer to this question is what is generally known as Birkhoff conjecture:  amongst all
convex billiards, the only {\it integrable} ones are the ones in ellipses (a circle is a distinct special case).\\

Despite its long history and the amount of attention that
this conjecture has captured, it remains essentially open.
As far as our understanding of integrable billiards is concerned,
the two most important related results are
the above--mentioned theorem by Bialy \cite{Bialy} (see also \cite{Woi}), a result by
Delshams and Ram\'irez-Ros \cite{DRR} in which they study 
entire perturbations of elliptic billiards and  prove that any nontrivial symmetric perturbation of the elliptic billiard is not integrable, and
a theorem by Mather \cite{Mather82} which proves the non-existence
of caustics (hence, the non-integrability) if the curvature of the boundary vanishes at one point.
This latter justifies the restriction of our attention  to strictly convex domains.\\

We shall see in the next subsection how this conjecture/question can be rephrased as a regularity question for Mather's $\beta$ function (see Question II bis).\\

\noindent {\bf 1.3 - Mather's minimal average action (or $\beta$-function) and billiards.}

At the beginning of the eighties Serge Aubry and John Mather developed, independently, what nowadays is commonly called {\it Aubry--Mather theory}.
This novel approach to the study of the dynamics of twist diffeomorphisms of the annulus, pointed out the existence of many {\it action-minimizing orbits}  for any given rotation number (for a more detailed introduction,  see  for example \cite{MatherForni, Siburg, SorLecNotes}).

More precisely, let $f: \R/\Z \times \R \longrightarrow \R/\Z \times\R$ a monotone twist map, {\it i.e.}, a $C^1$ diffeomorphism such that its lift  to the universal cover $\tilde{f}$ satisfies the following properties (we denote $(x_1,y_1)= \tilde{f}(x_0,y_0)$):
\begin{itemize}
\item[(i)] $\tilde{f}(x_0+1, y_0) = \tilde{f}(x_0, y_0) + (1,0)$,
\item[(ii)] $\frac{\partial x_1}{\partial y_0} >0$  (monotone twist condition),
\item[(iii)] $\tilde{f}$ admits a (periodic) generating function $h$ ({\it i.e.},  it is an exact symplectic map):
$$
y_1\,dx_1 - y_0\,dx_0 = dh(x_0,x_1).
$$
\end{itemize}
In particular, it follows from (iii) that:

\begin{equation} \label{genfuncttwistmap}
\left\{
\begin{array}{l}
y_1 = \frac{\partial h}{\partial x_1}(x_0,x_1)\\
y_0 = - \frac{\partial h}{\partial x_0}(x_0,x_1)\,.
\end{array}
\right.
\end{equation}

\begin{Rem}
The billiard map $f$ introduced above is an example of monotone twist map. In particular, its generating function (see (\ref{genfunctbill}))
is given by $h(x_0, x_1) = - \ell(x_0, x_1)$, where $\ell(x_0, x_1)$ denotes the euclidean distance between the two points on the boundary of the billiard domain corresponding
to $\gamma(x_0)$ and $\gamma(x_1)$.
\end{Rem}

As it follows from (\ref{genfuncttwistmap}), orbits $(x_i)_{i\in\Z}$ of the monotone twist diffeomorphism $f$ correspond to `critical points'  of the {\it action functional}
$$
\{x_i\}_{i\in\Z} \longmapsto \sum_{i\in \Z} h(x_i, x_{i+1}).
$$

Aubry-Mather theory is concerned with the study of orbits that minimize this action-functional amongst all configurations with a prescribed rotation number; recall that the rotation number of an orbit $\{x_i\}_{i\in\Z}$ is given by $\pi \omega = \lim_{i\rightarrow \pm \infty} \frac{x_i}{i}$, if this limit exists (in the billiard case, this definition leads to the same  notion of rotation number introduced in subsection 1.2). In this context, {\it minimizing} is meant in the statistical mechanical sense, {\it i.e.}, every finite segment of the orbit minimizes the action functional with fixed end-points.\\

\noindent {\bf Theorem (Aubry \& Mather).} {\it A monotone twist map possesses minimal orbits for every rotation number. For rational numbers there are always at least two periodic minimal orbits. Moreover, every minimal orbit lies on a Lipschitz graph over the $x$-axis.}\\

We can now introduce the {\it minimal average action} (or {\it Mather's $\beta$-function}).

\begin{Def}\label{defbeta}
Let $x^{\omega} = \{x_i\}_{i\in\Z}$ be any minimal orbit with rotation number $\omega$. Then, the value of the {\em minimal average action} at $\omega$ is given by (this value is well-defined, since it does not depend on the chosen orbit):
\begin{equation}\label{avaction}
\beta(\omega) = \lim_{N\rightarrow +\infty} \frac{1}{2N} \sum_{i=-N}^{N-1} h(x_i,x_{i+1}).
\end{equation}
\end{Def}

This function $\beta: \R \longrightarrow \R$ enjoys many properties and encodes interesting information on the dynamics. In particular:
\begin{itemize}
\item[i)] $\beta$ is strictly convex and, hence, continuous (see \cite{MatherForni});
\item[ii)] $\beta$ is differentiable at all irrationals (see \cite{Mather90});
\item[iii)] $\beta$ is differentiable at a rational $p/q$ if and only if there exists an invariant circle consisting of periodic minimal orbits of rotation number $p/q$ (see \cite{Mather90}).\\
\end{itemize}

In particular, being $\beta$ a convex function, one can consider its convex conjugate:
$$
\alpha( c ) = \sup_{\omega\in \R} \left[ \omega \, c - \beta(\omega)\right].
$$

This function -- which is generally called {\it Mather's $\alpha$-function} -- also plays an important r\^ole in the study minimal orbits and in Mather's theory (particularly in higher dimension, see for example  \cite{MS, SV}). We
refer interested readers to surveys \cite{MatherForni, Siburg, SorLecNotes}.\\

Observe that for each $\omega$ and $c$ one has:
$$
\alpha( c ) + \beta( \omega) \geq \omega  c,
$$
where equality is achieved if and only if $c\in \partial \beta(\omega)$ or, equivalently, if and only if  $\omega \in \partial \alpha (c )$  (the symbol $\partial$ denotes in this case the set of `subderivatives'  of the function, which is always non-empty and is a singleton if and only if the function is differentiable).\\

In the billiard case, since the generating function of the billiard map is the euclidean distance $-\ell$, the action of the orbit coincides -- up to a sign -- to the length of the trajectory that the ball traces on the  table $\Omega$. In particular, these two functions encode many dynamical  properties of the billiard (see \cite{Siburg} for more details):

\begin{itemize}
\item  For each $0<p/q\leq 1/2$, one has:
         $
         \beta(p/q) = - \frac{1}{q} {\mathcal ML}_\Om({p}/{q}).
         $
 \item $\beta$ is differentiable at $p/q$ if and only if there exists a caustic of rotation number $p/q$ ({\it i.e.}, all tangent orbits are periodic of rotation number $p/q$).
 \item If $\Gamma_{\omega}$ is a caustic with rotation number $\omega \in (0,1/2]$, then $\beta$ is differentiable at $\omega$ and $\beta'(\omega)= - {\rm length}(\Gamma_{\omega}) =: -|\Gamma_{\omega}|$ (see \cite[Theorem 3.2.10]{Siburg}). In particular,   $\beta$ is always differentiable at $0$ and $\beta'(0)= - |\partial \Omega|$.
 \item If $\Gamma_{\omega}$ is a caustic with rotation number $\omega \in (0,1/2]$, then one can associate to it another invariant, the so-called {\it Lazutkin invariant} $Q(\Gamma_{\omega})$.  More precisely 
\begin{equation}\label{lazinv}
Q(\Gamma_{\omega}) = |A-P| + |B-P| + |\stackrel \frown {AB} |
\end{equation}
 where $|\cdot |$ denotes the euclidean length and $|\stackrel \frown {AB}|$ the length of the arc on the caustic joining $A$ to $B$ (see figure \ref{lazutkin}).
 
This quantity is connected to the value of the $\alpha$-function. In fact, one can show that (see \cite[Theorem 3.2.10]{Siburg}):
 $$Q(\Gamma_\omega) = \alpha(  \beta'(\omega)  ) = \alpha( - |\Gamma_{\omega}|).$$
 
\begin{figure} [h!]
\begin{center}
\includegraphics[scale=0.8]{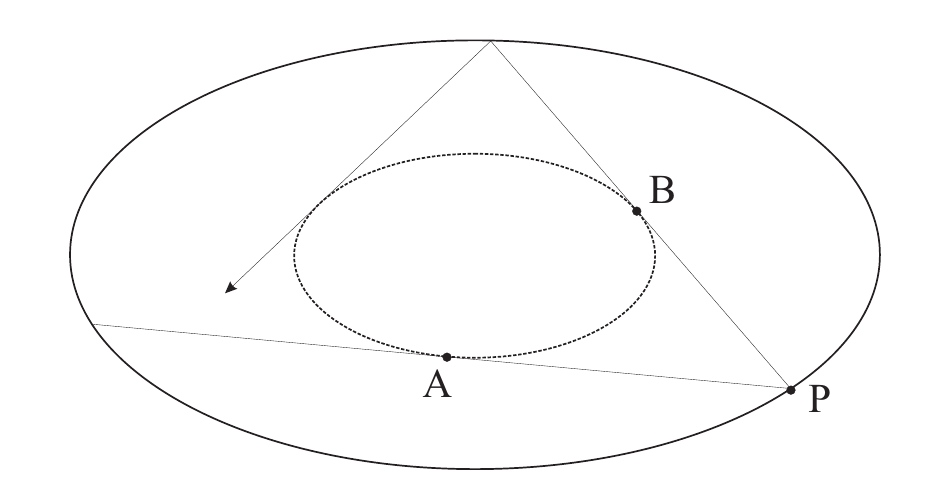}
\caption{Lazutkin invariant}
\label{lazutkin}
\end{center}
\end{figure}

 \end{itemize}

We can now rephrase Questions I and II (see above) in terms of these new objects.\\

\noindent{\bf Question I (bis).}  {\it Let $\Omega_1$ and $\Omega_2$ be two strictly convex planar domains with smooth boundaries and assume that 
$\beta_{\Omega_1} \equiv \beta_{\Omega_2}$. Is it true that $\Omega_1$ and $\Omega_2$ are isometric?}\\

Actually, one could ask even more. In fact, the knowledge of the dynamics near the boundary (for small angles) is sufficient to recover the curvature of the boundary and hence the global dynamics. Therefore:\\

\noindent{\bf Question I (ter).}  {\it Let $\Omega_1$ and $\Omega_2$ be two strictly convex planar domains with smooth boundaries and assume that 
$\beta_{\Omega_1}(\omega) =  \beta_{\Omega_2} (\omega)$ for all $\omega\in (0,\epsilon)$ for some small $\epsilon >0$. Is it true that $\Omega_1$ and $\Omega_2$ are isometric?}\\

\noindent{\bf Question II (bis).}  {\it Let $\Omega$ be a strictly convex planar domain with smooth boundary and assume that 
$\beta_{\Omega}$ is $C^\infty([0,\epsilon))$ for some small $\epsilon>0$.  Is it true that $\Omega$ is an ellipse?}\\

Observe that if $\beta_{\Omega}$ is $C^\infty([0,\epsilon))$, then  the billiard map is locally integrable near the boundary. In fact, $\beta$ will be differentiable at all rationals in $(0,\epsilon)$ and therefore there will be caustics corresponding to these rotation number. By semi-continuity arguments, one obtains caustics corresponding to irrational rotation number and hence a family of caustics that foliate a neighbourhood of the boundary.
Observe that if $\beta$ is differentiable in the whole domain of definition $(0,1/2]$, then it must be a circle by the aforementioned result by Bialy.
\\

\noindent {\bf 1.4 - Main results.}

Motivated by the above discussion, we would like to study more in depth the properties of Mather's $\beta$ and $\alpha$ functions and obtain explicit expressions for their (formal) Taylor expansions at, respectively, $\omega=0$ and $c=-\ell_0$ (where $\ell_0$ denotes the length of the boundary $\partial \Omega$). The coefficients in these expressions will be obtained only 
in terms of the curvature of the boundary (which, in fact, determines the dynamics univocally). 
The first order of these expressions have already appeared in \cite[Theorem 3.2.5]{Siburg}, but due to the nature of the argument (a perturbative argument), the analysis therein cannot be pushed further to higher orders. We shall follow here a different approach (more geometric), inspired by Amiran's work \cite{Amiran1}.\\

We shall prove the following.

\begin{Teo}\label{maintheorem}
Let $\Omega$ be a strictly convex planar domain with smooth boundary. Denote by $k(s)$ the curvature of $\partial \Omega$ with arc-length parametrization $s$. 
Let $\ell_0:= |\partial \Omega|$ be the length of the boundary and denote:
\begin{eqnarray*}
\I_1 &:=& \int_{0}^{\ell_0}ds  = \ell_0 \\
\I_3 &:=& \int_0^{\ell_0} k^{2/3} ds \\
\I_5 &:=& 
\int_0^{\ell_0}  \left( 9\ k^{4/3} +  \frac{8 \ \dot{k}^2}{k^{8/3}}  \right)ds \\
\I_7 &:=& 
 \int_0^{\ell_0} \left(
9\ k^2 + \frac{24\ \dot{k}^2}{k^2} + \frac{24\ \ddot{k}^2}{k^4}
- \frac{144\ \dot{k}^2 \ddot{k}}{k^5}
+ \frac{176\ \dot{k}^4}{k^6} 
\right) ds \\
\end{eqnarray*}
\begin{eqnarray*}
\I_9 &:=&
\int_0^{\ell_0}
\left[
\frac{281}{44800} k^{8/3}
+\frac{281\ \dot{k}^2}{8400\ k^{4/3}} +
\frac{167\ \ddot{k}^2}{4200\ k^{10/3}}
- \frac{167 \ \dot{k}^2 \ \ddot{k}}{700\ k^{13/3}}
+ \frac{{\dddot{k}}^2}{42\ k^{16/3}}
+\frac{559 \ \dot{k}^4}{2100\ k^{16/3}} \right.\\
&& \left. \;
-\; \frac{473 \ \ddot{k}^3}{4725\ k^{19/3}}
-\frac{10 \ \dddot{k}\ \dot{k} \ \ddot{k}}{21\ k^{19/3}}
+
\frac{5 \ \dddot{k} \ \dot{k}^3}{7\ k^{22/3}}
+ \frac{13142 \ \dot{k}^2\ \ddot{k}^2}{4725\ k^{22/3}}
-\frac{10777 \ \dot{k}^4 \ \ddot{k}}{1575\ k^{25/3}}
+\frac{521897 \ \dot{k}^6}{127575\
k^{28/3}}
 \right] ds.
\end{eqnarray*}

Then:
\begin{itemize}
\item
the formal Taylor expansion of $\beta$ at $\omega=0$,   $\beta(\omega) \sim \sum_{k=0}^{\infty}  \beta_k \frac{\omega^k}{k!}$, has coefficients: 
\begin{eqnarray*}
\beta_{2k}&=&0 \quad \mbox{for all $k$}\\
\beta_1 &=& - \I_1  \nonumber\\
\beta_3 &=&  \frac{1}{4}\ \I_3^3\nonumber \\
\beta_5 &=& - \frac{1}{144} \ \I_3^4\ \I_5 \\
\beta_7 &=& \frac{1}{320} \ \I_3^5 \left( \frac{14}{81} \I_5^2 - \I_3 \I_7  \right) \ = \
\frac{\I_3^5 \left(14\ 
   \I_5^2-81\ \I_3
   \I_7\right)}{25920}\nonumber\\
 \beta_9 &=& - 7\ \I_3^6 \left(  \I_3^2 \ 
   \I_9- \frac{1}{5600}
   \I_3\ \I_5 \
   \I_7+\frac{7 }{583200} \I_5^3 \right); \nonumber
   \end{eqnarray*}

\item the (formal) Taylor expansion of $(c+\ell_0)^{-3/2}\alpha( c )$ at $c=-\ell_0$ 
(note that $\alpha$ has in fact a square-root type singularity at the boundary),
${({c+\ell_0})^{-3/2}}{\alpha( c )}  \sim \sum_{k=0}^{\infty}  \alpha_k \frac{(c+\ell_0)^k}{k!}$, has coefficients:
\begin{eqnarray*}
\alpha_0 &=&   \frac{4\sqrt{2}}{3}  \I_3^{-3/2} \\
\alpha_1 &=& \frac{\sqrt{2}}{135} \I_3^{-7/2} \I_5\\
\alpha_2 &=& \frac{1}{56700 \sqrt{2}}\left(
\frac{72\, {\I_3} {\I_7}+7 \,{\I_5}^2}{
   {\I_3}^{11/2}}   \right)\\
\alpha_3 &=& \frac{1}{826686000 \sqrt{2}} \left(
 \frac{261273600\, {\I_3}^2 {\I_9}+ 21384\,  {\I_3}
   {\I_5} {\I_7}+1001\, {\I_5}^3}{
   {\I_3}^{15/2}}\right).\\
%
%
%
%
\end{eqnarray*}
\end{itemize}

\end{Teo}

\begin{Rem}
(1) The techniques used in the proof of the Theorem \ref{maintheorem}, allow one to obtain  explicit  expressions up to any arbitary high order (we restrict to order 11 just for the sake of this presentation).\\
(2) The coefficients $\beta_k$ are algebraically related to the set of spectral invariants introduced by Marvizi and Melrose \cite{MM} for strictly convex planar regions in order to investigate and give some partial answers to Kac's question on the isospectrality of planar domains. These computations provide explicit expressions for  those invariants as well (see the expressions for $\I_k$'s).
\end{Rem}

An easy consequence of these formulae is the following corollary.\\

\begin{Cor}\label{maincor}
Let $\Omega$ be a strictly convex planar domain with smooth boundary. Then:
$$
\beta_{3} + \pi^2 \beta_{1} \leq 0
$$
and equality holds if and only if $\Omega$ is a disc.\\
\end{Cor}

\begin{Rem}
In particular, the above corollary says that if the first two coefficients $\b_1$ and $\b_3$ coincide to those of the $\beta$-function of a disc, then the domain must be a disc. Therefore, the $\beta$-function univocally determines discs amongst all possible Birkhoff billiards.
It would be interesting to find a similar characterization for elliptic billiards. We can prove the following  result: the $\beta$-function  determines  univocally a given ellipse in the family of all ellipses.\\
\end{Rem}

\noindent {\bf Proposition  \ref{uniqueness}}. 
{\it If $\cE_1$ and $\cE_2$ are two ellipses such that $\beta_{\cE_1} \equiv \beta_{\cE_2}$, then $\cE_1$ and $\cE_2$ are the same ellipse. More generally: if the Taylor coefficients 
$ \beta_{\cE_1,1} =  \beta_{\cE_2,1}$ and $ \beta_{\cE_1,3} =  \beta_{\cE_2,3}$, then the same conclusion remains  true.}\\

The rest of the article is organized as follows. In Section \ref{secproof} we shall provide a proof of Theorem \ref{maintheorem}, which will be divided into several steps (subsections \ref{part1} -- \ref{part5}), while in subsection  \ref{part7}  Corollary \ref{maincor} will be deduced.
Finally, in Section \ref{examples} we shall discuss two families of examples: circular and elliptic billiards. In both case we shall provide expressions for Mather's $\beta$ functions and check the above formulae. In particular, in Section 3.2 we shall prove  Proposition \ref{uniqueness}.\\

\subsection*{Acknowledgements.} I would like to express my deepest gratitude to Vadim Kaloshin for having brought my attention to these (and many other) questions on billiards and for many interesting and engaging discussions. I wish to thank Corrado Falcolini for his precious help while using Mathematica for checking some of these computations.\\

\section{Proof of Theorem \ref{maintheorem}} \label{secproof}

In this section we prove Theorem \ref{maintheorem}.
Let $\Omega$ be a strictly convex  region in the plane bounded by a $C^{\infty}$ curve $\partial \Omega$, whose curvature is denoted by $k$ and whose radius of curvature by $\rho$.  
We aim at finding an expression of the (formal) Taylor expansion of $\beta$ at zero in terms of the curvature of the boundary.
In particular, if $\beta$ is smooth near $\omega=0$ (and consequently  the associated billiard map is integrable, \ie a neighbourhood of the boundary is smoothly foliated by caustics), this expansion will provide an expression of $\beta$ for sufficiently small rotation numbers.\\

 The proof will be splitted into several steps:

\begin{itemize}
\item[ $\S$\ref{part1} -] express the curvature of a caustic as a function of the  curvature of the boundary and the Lazutkin invariant;
\item[ $\S$\ref{part2} -] express the length of a  caustic as a function of the  curvature of the boundary and the Lazutkin invariant;
\item[ $\S$\ref{part3} -] express the rotation number of a caustic as a function of the  curvature of the boundary and its length;
\item[ $\S$\ref{part4} -] find  -- for rotation numbers for which a caustic exists -- an expression of $\beta$ as a function of the curvature of the boundary;
\item[ $\S$\ref{part5} -] discuss the existence of caustics near the boundary and find the (formal) Taylor expansion of $\beta$ at zero and other related quantities (for example, the $\alpha$ function, the relation between the rotation number and the Lazutkin invariant, etc $. . .$). End of the proof.\\
\end{itemize}

Moreover, in $\S$\ref{part7} we shall prove Corollary \ref{maintheorem}.\\

\subsection{Curvature of caustics and Lazutkin invariant} \label{part1} 
In this subsection we shall exploit some ideas already considered in \cite{Amiran1} and push them further to obtain information on the behaviour of higher order terms of the expansions (and correct some computational mistakes therein).

Let $\Gamma$ be a caustic and denote by $v$ its curvature, by $r$ its radius of curvature and by $L$ its Latzukin parameter. 
The first step consists in relating the curvature of $\partial \Omega$ to the curvature of $\Gamma$.

We identify smooth strictly invariant curves in $\R^2$ by their curvatures (see also \cite[Proposition 2.7]{MM}). To each closed curve we associate its curvature when the curve is parametrized by tangent angle (\ie the angle between tangent and $x$-axis), and to each positive $k\in C^{\infty}(\R/2\pi\Z,\R)$ with
$$
\int_0^{2\pi} k^{-1}(t)\cos t\,dt = \int_0^{2\pi} k^{-1}(t)\sin t\,dt =0,
$$
we associate the curve with coordinates
$$
x(\theta)= \int_0^{\theta} k^{-1}(t)\cos t\,dt  \quad {\rm and} \quad y(\theta)=\int_0^{\theta} k^{-1}(t)\sin t\,dt.
$$

Let us introduce the following parametrizations (we translate and rotate $\partial \Omega$ so that it passes through $(0,0)$ and its positive tangent direction at this point is $(1,0)$):

$$
\partial \Omega: \quad b(\varphi)= \left(\int_0^{\varphi} k^{-1}(t)\cos t\,dt, \; \int_0^{\varphi} k^{-1}(t)\sin t\,dt\right) \qquad \forall\; \varphi \in \R/2\pi\Z.
$$
and
$$
\Gamma:  \quad  a(\theta) = \left(x^0_{\G}+\int_0^{\theta} v^{-1}(t)\cos t\,dt, \; y^0_\G+ \int_0^{\theta} v^{-1}(t)\sin t\,dt\right) \qquad \forall\;\theta \in \R/2\pi\Z.
$$

Since $\G$ is a caustic of $\partial \Omega$, we can say that $\partial \Omega$ is an {\it $L$-evolute} of $\Gamma$, where $L=L(\Gamma)$ is the Latzukin parameter of the caustic $\G$ (see definition \ref{lazinv} subsection 1.3). Therefore, for each $\f \in \R/2\pi\Z$ there exist $\th_1,\th_2 \in \R/2\pi\Z$ and $t_1,t_2 >0$ (see figure \ref{figureamiran}) such that:
\beqa{eqb}
&& b(\f) = a(\th_1) + t_1(\cos \th_1, \sin \th_1) = a(\th_2) - t_2(\cos \th_2, \sin \th_2) \quad  {\rm and}\\
&& L=t_1+t_2 - (s(\th_2)-s(\th_1)),\nonumber
\eeqa
where  $s(\th):=\int_0^\th v^{-1}(t)$  denotes the arc-length along $a$ between $a(0)$ and $a(\th)$:

\begin{figure} [h!]
\begin{center}
\includegraphics[scale=0.7]{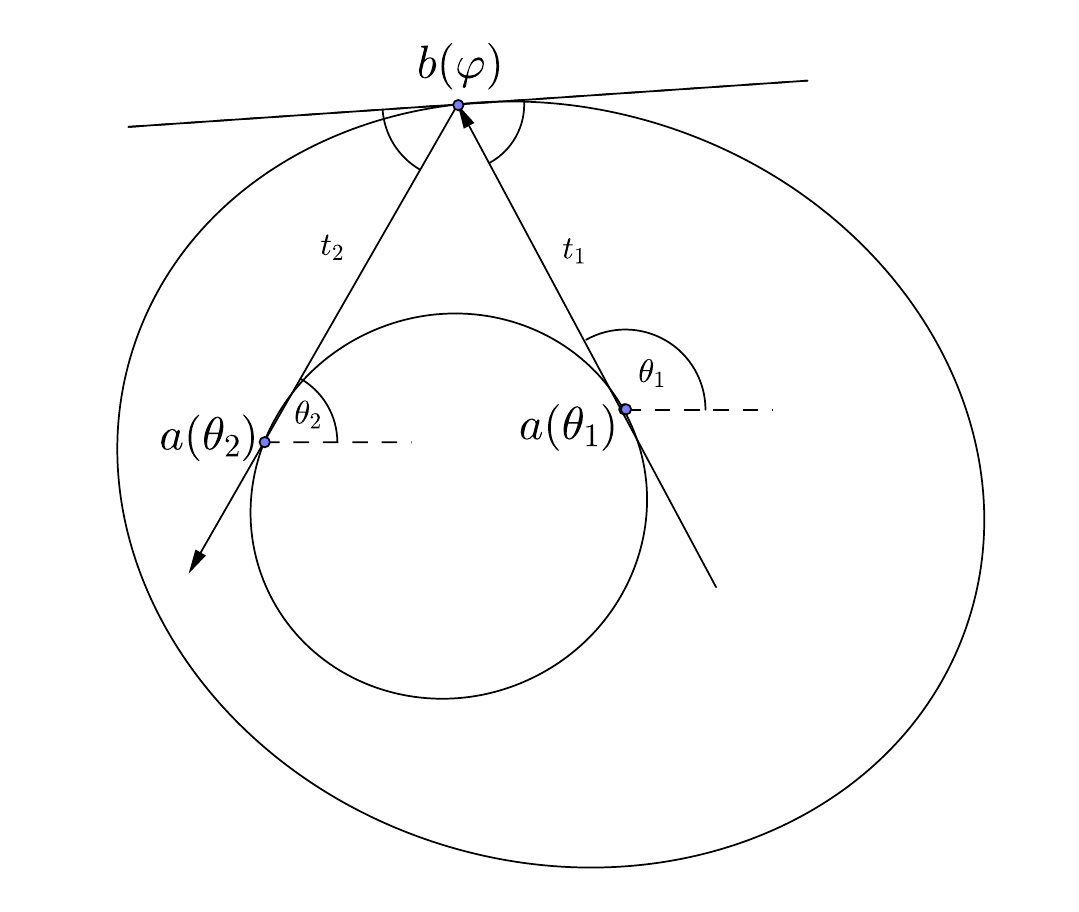}
\caption{}
\label{figureamiran}
\end{center}
\end{figure}

Since $\G$ is assumed to be an invariant curve for the billiard map on $\partial \Omega$, then one can deduce that $\f = \frac{\th_1+ \th_2}{2}$ (see figure). Moreover:
\beqa{relationt1t2}
\left\{
\begin{array}{l}
t_1 \cos \th_1 + t_2 \cos \th_2 = \int_{\th_1}^{\th_2} v^{-1}(t)\cos t\, dt\\
t_1 \sin \th_1 + t_2 \sin \th_2 = \int_{\th_1}^{\th_2} v^{-1}(t)\sin t\, dt.
\end{array}
\right.
\eeqa
It follows from above that:\footnote{It is sufficient to expand $(t_1 \cos \th_1 + t_2 \cos \th_2)(\sin \th_1 - \sin \th_2) + (t_1 \sin \th_1 + t_2 \sin \th_2)(\cos \th_2 - \cos \th_1)$ and simplify.}
\beqano
t_1+t_2 &=& \frac{1}{\sin (\th_2-\th_1)} \int_{\th_1}^{\th_2} v^{-1}(t)[\sin \th_2 \cos t - \cos \th_2 \sin t + \cos \th_1 \sin t - \sin \th_1 \cos t]\,dt\, \\
&=& \frac{1}{\cos \D} \int_{\f-\D}^{\f+\D} \cos (\f - t)v^{-1}(t)\,dt,
\eeqano
where 
\beqa{defDelta}
\D:= \th_2-\f = \f-\th_1 \qquad \mbox{and therefore}\qquad \th_2-\th_1=2 \Delta.
\eeqa

In particular, this shows that (use the change of variable $t=u+\f$):
{\small
\beqa{eqQ}
L &=& t_1+t_2 - (s(\th_2)-s(\th_1)) = t_1+t_2 - \int_{\th_1}^{\th_2} v^{-1}(t)\,dt  \\
&=&  \frac{1}{\cos \D} \int_{\f-\D}^{\f+\D} \cos (\f - t)v^{-1}(t)\,dt - \int_{\th_1}^{\th_2} v^{-1}(t)\,dt  \nonumber\\
&=& \frac{1}{\cos \D} \int_{-\D}^{\D} \cos(u) \, v^{-1}(\f+u)\,du - \int_{-\D}^{\D} v^{-1}(\f+u)\,du\,,  \nonumber\\
&=&  \frac{1}{\cos \D} \; \int_{0}^{\D} \cos(u)\,\left(v^{-1}(\f+u) + v^{-1}(\f-u)\right)\,du -
\int_{0}^{\D} \left(v^{-1}(\f+u) + v^{-1}(\f-u)\right)\,du.\nonumber
\eeqa}





Expanding in $\D$, we obtain:
\beqa{eqQ2}
L &=& 
\frac{2}{3} r(\f)\, \D^3 +\frac{1}{15} 
   \left[ r''(\f)+4 \, r(\f)\right]\, \D^5 + \left[\frac{  3\, r^{(4)}(\f)+32\, r''(\f)+136\,  r(\f)}{1260} \right]  \D^7  \nonumber\\
   &&\; +\;       \left[\frac{
r^{(6)}(\f) + 20\, r^{(4)}(\f) + 232\,  r''(\f)+992\,   r(\f)}{22680} \right] \D^9 + O\left(\D^{11}\right  ).
\eeqa

We can now invert the above expression and obtain an expansion of $\Delta$ in terms of $L$ (we write $r$ instead of $r(\f)$):
{\small
\beqa{eqdelta}
&& \D \; =\;  
\left(\frac{3}{2}\right)^{1/3} r^{-1/3} L^{1/3}
 +
   \left[ \frac{
   -r''-4 r}{20
   r^2} \right] L  \\
&& \;+ \;   
\left(\frac{3}{2}\right)^{2/3}
   \left[ \frac{
   -15 \ r^{(4)}\, r+288\ r\
   r''+56 \ r''^2 + 216\
   r^2} {8400\  r^{11/3} } \right]  L^{5/3}  \nonumber\\ 
&&\; + \;   %
\left(\frac{3}{2}\right)^{1/3}
\frac{
   \left(-5 \ r^{(6)} \ r^2 + 260\
   r^{(4)} \ r^2-1976\ r^2\
   r'' - 1224\ r \ r''^2 - 182\
   r''^3 + 90\ r^{(4)}\ r\ r'' - 288\
   r^3\right)}{100800\
r^{16/3}
   }\ L^{7/3}   \nonumber \\
   && \;+\;  O\left(L^3\right).\nonumber
\eeqa
}

The curvature of $\partial \Omega$ at a point $b(\f)=(x(\f), y(\f))$ is given by:
$$
k(\f) = \left(\left(\frac{dx}{d\f}\right)^2 + \left(\frac{dy}{d\f}\right)^2   \right)^{-\frac{1}{2}}. 
$$
In particular, it follows from (\ref{eqb}) and the definition of $\D=\f-\th_1=\th_2-\f$ that:
$$
\left\{\begin{array}{l}
\frac{dx}{d\f} = \frac{\cos (\f-\D)}{v(\f-\D)} \cdot \frac{d(\f-\D)}{d\f} - t_1 \sin (\f-\D) \cdot \frac{d(\f-\D)}{d\f} + \cos(\f-\D)\cdot \frac{d t_1}{d\f}\\
\frac{dy}{d\f} = \frac{\sin (\f-\D)}{v(\f-\D)} \cdot  \frac{d(\f-\D)}{d\f} + t_1 \cos (\f-\D)\cdot\frac{d(\f-\D)}{d\f} + \sin(\f-\D)\cdot \frac{d t_1}{d\f}.
\end{array}\right.
$$

Therefore, 
\beqa{formulacurvature}
\left(\frac{dx}{d\f}\right)^2 + \left(\frac{dy}{d\f}\right)^2  = \left[
v^{-1}(\f-\D) \left(1- \frac{d\D}{d\f} \right) + \frac{dt_1}{d\f} \right]^2 + t_1^2 \left(1- \frac{d\D}{d\f} \right)^2.
\eeqa

Let us express \footnote{Observe that the corresponding formula in \cite[p.352]{Amiran1} is not correct due to some computational mistake.} this quantity in terms of $L$.

First of all, it follows from (\ref{relationt1t2}) that:

{\small
\beqa{t1bis}
t_1 &=& - \frac{1}{2\sin \D} \int_{-\D}^{\D} \sin (u)\ v^{-1}(u+\f)\,du  + 
\frac{1}{2\cos \D} \int_{-\D}^{\D} \cos (u)\ v^{-1}(u+\f)\,du   \nonumber\\
&=& -\; \frac{1}{2\sin \D} \int_{0}^{\D} \sin (u) \left[v^{-1}(u+\f) - v^{-1}(u-\f)\right]\,du
 \nonumber\\
&& \; +\;
\frac{1}{2\cos \D} \int_{0}^{\D} \cos (u) \left[v^{-1}(u+\f) + v^{-1}(u-\f)\right]\,du  \nonumber\\
&=& r\ \D - \frac{1}{3} r'\ \D^2 + \frac{1}{6}\left[ r'' +2 r \right]\ \D^3
   + \frac{1}{90}  \left[-3
   r^{(3)}-2 r'\right]\ \D^4 +
     +\frac{1}{120}
   \left[ r^{(4)}+4 r''+16
   r\right]\ \D^5 
   + \nonumber\\
   && \;+\;
    \left[ \frac{-9
   r^{(5)}-12 r^{(3)}-16
   r'}{7560}\right] \D^6
    + \left[\frac{
   r^{(6)}+6 r^{(4)}+64
   r''+272 r}{5040}\right] \D^7 \nonumber\\
   &&\;+\;
   \left[\frac{
   -32 r^{(3)}-48 r'-5
   \left(r^{(7)}+2
   r^{(5)}\right)}{226800} \right] \D^8 + O\left(\D
   ^{9}\right).
\eeqa
}

Recalling (\ref{eqdelta}), we also obtain (we write $r$ instead of $r(\f)$):
{\small
\beqa{vmenodelta}
&& v^{-1}(\f-\D) \ =\  r(\f-\D)\ =\ 
r - \left(\frac{3}{2}\right)^{1/3}
\left[
\frac{r'}{r^{1/3}}\right] L^{1/3}
   +
      \left(\frac{3}{2}\right)^{2/3}  \left[
      \frac{r''}{2   r^{2/3}}\right]   L^{2/3}
  \nonumber\\ 
&&\; +\;   \left[ \frac{r'
   \left(r''+4 r\right)-5 r
   r^{(3)}}{20 \ r^2} \right] 
   L
   +
     \left(\frac{3}{2}\right)^{1/3}
   \left[
   \frac{ 5 r r^{(4)}-4
   r'' \left(r'' +4
   r\right)}{80
   r^{7/3}}  \right]     L^{4/3} \nonumber\\
   &&\; +\;
   \left(\frac{3}{2}\right)^{2/3} \left[  \frac{ -56
   r' r''^2-3 r^2 \left(72
   r'+35 \left(r^{(5)}-8
   r^{(3)}\right)\right)+3 r
   \left(70 r^{(3)} r''+r'
   \left(5 r^{(4)}-96
   r''\right)\right)}{8400\
   r^{11/3}} \right]
      L^{5/3}\nonumber\\
   &&\;+\; \left[\frac{126
   r''^3+2 r \left(344 r''-85
   r^{(4)}\right) r''+r^2
   \left(656 r''+35 \left(r^{(6)}-16
   r^{(4)}\right)\right)}{11200
   r^4} \right]
   L^2 
   \nonumber\\
   &&\;+\;
   \left(\frac{3}{2}\right)^{1/3} 
   \left[ 
   \frac{ 182 r'
   r''^3+9 r^3 \left(32 r'-5
   \left(r^{(7)}-28 r^{(5)}+88
   r^{(3)}\right)\right)+18 r r''
   \left(r' \left(68 r''-5
   r^{(4)}\right)-35 r^{(3)}
   r''\right)}{100800 \ 
   r^{16/3}}  \right.\nonumber\\
&&\qquad \qquad \quad +\; \left.   \frac{
    r^2 \left(r'
   \left(5 r^{(6)}-260 r^{(4)}+1976
   r''\right)+45 \left(7 r^{(5)}
   r''+r^{(3)} \left(3 r^{(4)}-80
   r''\right)\right)\right)}{100800 \
   r^{16/3}}\right]
   L^{7/3}\nonumber\\
   && \; - \;  \left(\frac{3}{2}\right)^{2/3} \left[
   \frac{8624
   r''^4+528 r \left(116 r''-25
   r^{(4)}\right) r''^2+8 r^2
   \left(225 {r^{(4)}}^2+13984
   {r''}^2+\left(340 r^{(6)}-8320
   r^{(4)}\right) r''\right)}
   {4032000\
   r^{17/3}}
   \right.\nonumber\\
   && \qquad \qquad \quad \left. 
   +\; \frac{
   9
   r^3 \left(-25 r^{(8)}+1120
   r^{(6)}-7360 r^{(4)}+3584
   r''\right)}{4032000\
      r^{17/3}} \right]
   L^{8/3}
   +O\left(L^3\right).
         \eeqa
}

Moreover, it follows from (\ref{eqdelta}), (\ref{t1bis}), and the fact that $L$ is constant with respect to $\f$ (since $\G$ is a caustic) that
\beqa{derivativeDelta}
\frac{d\D}{d\f} &=& 
-  \left(\frac{3}{2} \right)^{1/3}  \left[\frac{
   r'}{ 3\ r^{4/3}}
   \right]
   L^{1/3}
    +
   \left[\frac{-r^{(3)}-4 r'}{20
   r^2}-\frac{r' \left(-r''-4
   r\right)}{10\  r^3}\right]
   L\\
   &&\;+\
  \left(\frac{3}{2} \right)^{2/3} \left[ 
   \frac{-15 r r^{(5)}+288
   r r^{(3)}+432 r r'-15
   r^{(4)} r'+112 r^{(3)}
   r''+288 r' r''}{8400\
 r^{11/3}} \right. \nonumber\\
&& \qquad \qquad \quad \left. - \;
   \frac{11
   r' \left(-15 r^{(4)} r+288
   r r''+56 r''^2+216
   r^2\right)}{25200\ r^{14/3}} \right]
          L^{5/3} \nonumber\\
    && \; + \; \left(\frac{3}{2} \right)^{1/3}
   \left[
   \frac{ -5 r^{(7)} r^2+260
   r^{(5)} r^2-1976 r^{(3)}
   r^2-864 r^2 r'-10 r^{(6)}
   r r'+90 r^{(5)} r r''+ 90
   r^{(3)} r^{(4)} r}{100800  \ r^{16/3}} \right. \nonumber\\
   &&\qquad \qquad \quad \left.
   +\; \frac{520
   r^{(4)} r r'-2448 r^{(3)}
   r r''-546 r^{(3)}
   r''^2-3952 r r' r''-1224
   r' r''^2+90 r^{(4)} r'
   r''}
   {100800  \ r^{16/3}} \right. \nonumber\\
   &&\qquad\qquad\quad
   \left.
   -
   \frac{ r' \left(-5 r^{(6)}
   r^2+260 r^{(4)} r^2-1976
   r^2 r''-1224 r r''^2-182
   r''^3+90 r^{(4)} r r''-288
   r^3\right)}{18900 \ r^{19/3}}
   \right]
   L^{7/3} \nonumber\\
&& \;+\;   O\left(L^3\right) \nonumber
\eeqa
and
{\small
\beqa{derivativet1}
\frac{dt_1}{d\f} &=&
\left(\frac{3}{2}\right)^{-2/3}
\left[ \frac{    r'}{r^{1/3}}
   \right] L^{1/3}
   +
\left(\frac{3}{2}\right)^{2/3} 
\left[   \frac{
   2 r'^2-3 r
   r''}{9\ 
   r^{5/3}} \right]
   L^{2/3}
   + \left[\frac{ r
   r^{(3)}-r' r''}{5 \ 
   r^2}\right] L \\
   &&
   +\;
   \left(\frac{3}{2}\right)^{1/3} 
   \left[
   \frac{ 9 r^2
   \left(2 r''-r^{(4)}\right)-14
   r'^2 r''+6 r \left(r''^2-4
   r'^2+3 r^{(3)}
   r'\right)}{180 \} r(f)^{10/3}}
   \right]
   L^{4/3} \nonumber\\
   &&+\;
   \left(\frac{3}{2}\right)^{2/3} 
   \left[ \frac{
   616 r' r''^2-9 r^2
   \left(-15 r^{(5)}+92 r^{(3)}+24
   r'\right)+3 r \left(5 r'
   \left(92 r''-15
   r^{(4)}\right)-154 r^{(3)}
   r''\right)}{12600\  r^{11/3}}
   \right]
   L^{5/3} \nonumber\\
   &&
   + \; \left[
   \frac{
    168 r'^2 r''^2-3 r^3
   \left(5 r^{(6)}-68 r^{(4)}+32
   r''\right)+6 r \left(-7
   r''^3+r'^2 \left(96 r''-5
   r^{(4)}\right)-42 r^{(3)} r'
   r''\right)}{5600 \ 
   r^5} \right. \nonumber\\
   &&\qquad\quad + \left. \frac{
   2 r^2 \left(28
   {r^{(3)}}^2-96 r''^2+96 r'^2+33
   r^{(4)} r''+20
   \left(r^{(5)}-15 r^{(3)}\right)
   r'\right)}{5600 \ 
   r^5} 
   \right] L^2 \nonumber\\
   &&+\;
   \left(\frac{3}{2}\right)^{1/3} \left[
   \frac{-728 r'
   r''^3+3 r^3 \left(5
   r^{(7)}-125 r^{(5)}+248
   r^{(3)}+48 r'\right)+18 r
   r'' \left(28 r^{(3)} r''+5
   r' \left(5 r^{(4)}-33
   r''\right)\right)}{37800\   r^{16/3}} \right. \nonumber\\
   && \qquad \qquad \quad \left. + \; \frac{
   -r^2 \left(7
   r' \left(5 r^{(6)}-125
   r^{(4)}+248 r''\right)+27 \left(5
   r^{(5)} r''+r^{(3)} \left(5
   r^{(4)}-66
   r''\right)\right)\right)}{37800\   r^{16/3}}
   \right]
   L^{7/3} \ + \nonumber\\
   \end{eqnarray}
   }
{\small   \begin{eqnarray}
&&   +\;  \left(\frac{3}{2}\right)^{2/3}
\left[
\frac{ -36652
   r'^2 r''^3-9 r^4 \left(25
   r^{(8)}-1030 r^{(6)}+5344
   r^{(4)}+48 r''\right)}{4536000\    r^{20/3}}
   \right.\nonumber\\
 &&
 \qquad\qquad\quad + \left.
   \frac{924 r
   r'' \left(7 r''^3+3 r'^2
   \left(5 r^{(4)}-68 r'' \right)+ 63
   r^{(3)} r' r'' \right)-22
   r^2 \left(27 r'' \left(28
   {r^{(3)}}^2-68 r''^2+19 r^{(4)}
   r''\right)\right)}{4536000\    r^{20/3}}
   \right.\nonumber\\
 &&
 \qquad\qquad\quad + \left.
   \frac{-22
   r^2 \left(
   r'^2 \left(9232
   r''+25 \left(r^{(6)}-52
   r^{(4)}\right)\right)+135 r'
   \left(3 r^{(3)} r^{(4)}+4
   \left(r^{(5)}-21 r^{(3)}\right)
   r''\right)\right)
   }{4536000\    r^{20/3}}
   \right.\nonumber\\
 &&
 \qquad\qquad\quad + \left.
   \frac{
   6 r^3 \left(270
   {r^{(4)}}^2+9232 r''^2+192
   r'^2+27 r^{(3)} \left(25
   r^{(5)}-284
   r^{(3)}\right)+\left(430
   r^{(6)}-8968 r^{(4)}\right)
   r''\right)
    }{4536000\    r^{20/3}}
   \right.\nonumber\\
 &&
 \qquad\qquad\quad + \left.
   \frac{
   6 r^3 \left(
   \left(125 r^{(7)}-5420
   r^{(5)}+30608 r^{(3)}\right)
   r'\right)}{4536000\    r^{20/3}}
   \right]
   L^{8/3} \nonumber\\
   && +\;  O\left(L^3\right). \nonumber
       \eeqa
 }

%

Let us now substitute these estimates in (\ref{formulacurvature}) and consider its Taylor expansion:
{\small
\beqa{expressionradiuscurvature}
\rho(\f) &=&  k^{-1}(\f) \;=\;  \sqrt{\left[
r(\f-\D) \left(1- \frac{d\D}{d\f} \right) + \frac{dt_1}{d\f} \right]^2 + t_1^2 \left(1- \frac{d\D}{d\f} \right)^2} = \\
&=& r(\f) + C_1[r(\f)] \cdot L^{2/3} + C_2[r(\f)] \cdot  L^{4/3} + C_3[r(\f)] \cdot  L^2 + C_4[r(\f)] \cdot  L^{8/3} + O\left(L^{10/3}\right) \nonumber
\eeqa
}
where $C_1, C_2, C_3, C_4: C^{\infty}(\sfrac{\R}{2\pi\Z}, \R_+) \longrightarrow C^{\infty}(\sfrac{\R}{2\pi\Z}, \R)$ are operators given by:

{\small
\begin{eqnarray*}
C_1[r(\f)] &:=& \left(\frac{3}{2}\right)^{2/3}
\left[\frac{3 r \left(r'' +3 r \right)-2 r'^2}{18\ r^{5/3}} \right] \\
C_2[r(\f)] &:=& \left(\frac{3}{2}\right)^{1/3}
\left[ \frac{9 r^2 \left(r^{(4)}-2
   r''\right)+28 r'^2 r''+12 r
   \left(-r''^2+3 r'^2-2 r^{(3)}
   r'\right)+81 r^3}{720\ 
r^{10/3}}\right]
\\
C_3[r(\f)] &:=&
-\frac{168 r'^2 r''^2+r^3
   \left(-5 r^{(6)}+77 r^{(4)}+73
   r''\right)+6 r \left(-7
   r''^3+r'^2 \left(82 r''-5
   r^{(4)}\right)-28 r^{(3)} r'
   r''\right)}{11200 \ r^5}
   \\
   && \; -\;
   \frac{
    2 r^2 \left(14
   {r^{(3)}}^2-75 r''^2+9 r'^2+19
   r^{(4)} r''+2 \left(5
   r^{(5)}-82 r^{(3)}\right)
   r'\right)+9 r^4}{11200 \ r^5}
\\
%
C_4[r(\f)] &:=&
\left(\frac{3}{2}\right)^{2/3}
\left[
\frac{146608\  r'^2 r''^3+9 r^4
   \left(25 r^{(8)}-1140
   r^{(6)}+4638 r^{(4)}+12988
   r''\right)}{36288000\
   \ r^{20/3}}
   \right. \\
   && \qquad \qquad \quad \left. +\; \frac{
   - 3696 \ r r'' \left(7
   r''^3+3 r'^2 \left(5
   r^{(4)}-61 r''\right)+42
   r^{(3)} r' r''\right)}{36288000\
   \ r^{20/3}}
   \right. \\
   && \qquad \qquad \quad \left. +\; \frac{
    88
   r^2 \left(18 r'' \left(21
   {r^{(3)}}^2-88 r''^2+18 r^{(4)}
   r''\right)\right)
   }{36288000\
   \ r^{20/3}}
   \right. \\
   && \qquad \qquad \quad \left. +\; \frac{
    88
   r^2 \left(
   r'^2 \left(25
   r^{(6)}-1165 r^{(4)}+5803
   r''\right)+54 r' \left(5
   r^{(3)} r^{(4)}+\left(5
   r^{(5)}-122 r^{(3)}\right)
   r''\right)\right)}{36288000\
   \ r^{20/3}}
   \right. \\
   && \qquad \qquad \quad \left. -\; \frac{
   24 r^3
   \left(135 {r^{(4)}}^2+4156
   r''^2+7185 r'^2+54 r^{(3)}
   \left(5 r^{(5)}-61
   r^{(3)}\right)+\left(160
   r^{(6)}-4324 r^{(4)}\right)
   r''\right)}{36288000\
   \ r^{20/3}}
   \right. \\
   && \qquad \qquad \quad \left. -\; \frac{
   24 r^3 \left(
   2 \left(25 r^{(7)}-1165
   r^{(5)}+5803 r^{(3)}\right)
   r'\right)+ 38799\ r^5}{36288000\
   \ r^{20/3}}
   \right].
\end{eqnarray*}
}

Next goal is to invert the above expression and write $r(\f)$ in terms of $\rho(\f)$. 
{\small
\beqa{eqk2}
r(\f) &=& \rho(\f) - C_1[r(\f)] \cdot L^{2/3} - C_2[r(\f)] \cdot  L^{4/3} - C_3[r(\f)] \cdot  L^2 - C_4[r(\f)] \cdot  L^{8/3} + O\left(L^{10/3}\right) \nonumber \\
&=& \rho(\f) + 
 C_1\left[ \rho(\f) - C_1\left[ \rho(\f) - C_1\left[  \rho(\f) - C_1\left[\rho(\f)\right] \cdot L^{2/3}           \right] \cdot L^{2/3} - C_2\left[\rho(\f)\right] \cdot  L^{4/3}\right] \cdot L^{2/3} \right. \nonumber \\
 &&\qquad \qquad \quad \left. 
 - C_2\left[\rho(\f) - C_1\left[\rho(\f)\right] \cdot L^{2/3}\right] \cdot  L^{4/3} - C_3\left[\rho(\f)\right] \cdot  L^2 \right] \cdot L^{2/3}  \nonumber\\
&& \; +\;  C_2\left[\rho(\f) - C_1\left[\rho(\f) - C_1\left[\rho(\f)\right] \cdot L^{2/3}\right] \cdot L^{2/3} - C_2\left[\rho(\f)\right] \cdot  L^{4/3} \right] \cdot  L^{4/3}  \nonumber\\
&& \; +\; C_3\left[\rho(\f) - C_1[\rho(\f)\right ] \cdot L^{2/3}] \cdot  L^2 \nonumber\\
&& \; +\; C_4\left[\rho(\f)\right] \cdot  L^{8/3} + O\left(L^{10/3}\right) \nonumber\\
&=& \rho(\f) - A(\f)\cdot L^{2/3} + B(\f)\cdot L^{4/3} + C(\f)\cdot L^{2} + D(\f)\cdot L^{8/3} + O\left(L^{10/3}\right)
\eeqa
}
where:
\begin{eqnarray*}
A(\f) &:=&   - \left(\frac{3}{2}\right)^{2/3} \left[\frac{2 \rho'^2-3 \rho \left(\rho''+3 \rho\right)}{18 \ \rho^{5/3}}\right] \\
&=&  \left(\frac{3}{2}\right)^{2/3} \left[\frac{1}{2} \rho^{1/3} + \frac{1}{2}\frac{d^2}{d\f^2}\left(\rho^{1/3}\right)\right]\\
&=&  \left(\frac{3}{2}\right)^{2/3} \frac{1}{2} \rho^{1/3} + \frac{d f_A}{d\f}\\
B(\f) &:=&  \left( \frac{3}{2}\right)^{1/3} \frac{1}{720}
\left[\frac{9}{\rho^{1/3}} +  8 \frac{\rho'^2}{\rho^{7/3}} \right]  + \frac{df_B}{d\f}\\
C(\f) &:=& \frac{1}{11200} \left[
\frac{9}{\rho} + 24 \left(\frac{\rho'^2 + \rho''^2 }{\rho^3}\right)
- 40 \frac{ \rho'^4}{\rho^5}\right] +  \frac{df_C}{d\f}\\
D(\f) &:=& \left( \frac{3}{2}\right)^{2/3} \frac{1}{90} 
\left[
\frac{281}{44800}\cdot \frac{1}{\rho^{5/3}} +
\frac{1}{\rho^{11/3}} \left(
\frac{281}{8400}\rho'^2
+ \frac{167}{4200} \rho''^2
+
\frac{{\rho^{(3)}}^2}{42}
\right) 
+\frac{473 }{4725}\cdot
\frac{\rho''^3}{\rho^{14/3}}
\right.\\
&&\qquad \qquad \quad \left.
- \frac{1}{\rho^{17/3}} \left(
\frac{11}{120} \rho'^4 + 
\frac{473}{945} \rho'^2 \rho''^2
\right)
+ \frac{781}{1458} \frac{\rho'^6}{\rho^{23/3}}
\right] + \frac{df_D}{d\f}
 \end{eqnarray*} 
and $f_A, f_B, f_C, f_D \in C^{\infty}(\sfrac{\R}{2\pi \Z})$ given by:
{\small
\begin{eqnarray*}
f_A (\f) &=& \frac{1}{2} \left(\frac{3}{2}\right)^{2/3} \frac{d}{d\f}\left( \rho^{1/3}\right) = \frac{1}{6}\left(\frac{3}{2}\right)^{2/3}  \frac{\rho'}{\rho^{2/3}}\\
f_B(\f) &=& \frac{1}{2160} \left( \frac{3}{2}\right)^{1/3}  \left[ \frac{140 \rho'^3+9 \rho^2
   \left(7 \rho^{(3)} +6 \rho'\right)-204 \rho\rho'\rho''}{\rho^{10/3}}\right]%
 \\
f_C(\f) &=& 
-\frac{31 \rho^{(5)}}{6720 \rho^2}-\frac{7 \rho^{(3)}}{4800
   \rho^2}-\frac{25 \rho'^5}{108
   \rho^6}+\frac{\rho'^3}{30240
   \rho^4}+\frac{79 \rho'}{33600
   \rho^2}+\frac{23 \rho^{(4)}
   \rho'}{672 \rho^3}+\frac{73
   \rho^{(3)} \rho''}{1200 \rho^3}-\frac{23 \rho^{(3)} \rho'^2}{144 \rho^4}\\
   && \;+\; \frac{283 \rho'^3 \rho''}{540 \rho^5}-\frac{1607 \rho' \rho''^2}{7200 \rho^4}+\frac{19 \rho' \rho''}{16800 \rho^3}\\
f_D(\f) &=&  \frac{1}{3} \left( \frac{3}{2}\right)^{2/3} \left[ 
\frac{127 \rho^{(7)}}{89600 \rho^{8/3}}-\frac{31 \rho^{(5)}}{67200 \rho^{8/3}}-\frac{67 \rho^{(3)}}{96000 \rho^{8/3}}+\frac{211945 \rho'^7}{52488 \rho^{26/3}}-\frac{9613 \rho'^5}{97200 \rho^{20/3}}-\frac{\rho'^3}{37800
   \rho^{14/3}}+\frac{839 \rho'}{1008000 \rho^{8/3}} \right.\\
   &&  \left.\;-\;\frac{207
   \rho^{(6)} \rho'}{11200 \rho^{11/3}}-\frac{4661 \rho^{(5)}
   \rho''}{100800 \rho^{11/3}}+\frac{44473 \rho^{(5)}
\rho'^2}{302400 \rho^{14/3}}-\frac{41 \rho^{(3)}\rho^{(4)}}{576 \rho^{11/3}}-\frac{389257 \rho ^{(4)}
   \rho'^3}{453600 \rho^{17/3}}+\frac{89 \rho^{(4)}
   \rho'}{12600 \rho^{11/3}}+\frac{121393 \rho^{(3)}
   \rho''^2}{216000 \rho^{14/3}} \right.\\
&&\left.\;   +\;\frac{197 \rho^{(3)}
   \rho''}{18000 \rho^{11/3}}+\frac{1034933 \rho^{(3)} \rho '^4}{272160 \rho^{20/3}}-\frac{31 \rho ^{(3)}
   \rho'^2}{672 \rho^{14/3}}+\frac{6763 {\rho^{(3)}}^2
   \rho'}{16800 \rho^{14/3}}-\frac{56771 \rho'^5
   \rho''}{4536 \rho^{23/3}}+\frac{14493167 \rho'^3
   \rho''^2}{1360800 \rho^{20/3}}+\frac{1621 \rho'^3 \rho''}{8400 \rho^{17/3}}\right.\\
&&\left.\;   -\; \frac{14553127 \rho'
   \rho''^3}{6804000 \rho^{17/3}}-\frac{3421 \rho' \rho''^2}{54000 \rho^{14/3}}+\frac{3
   \rho' \rho''}{4000 \rho^{11/3}}+\frac{63689 \rho^{(4)}
   \rho' \rho''}{100800 \rho^{14/3}}-\frac{1040447 \rho^{(3)} \rho'^2 \rho''}{226800 \rho^{17/3}} \right].
\end{eqnarray*}
}


\subsection{Length of a  caustic as a function of the  curvature of the boundary and the Lazutkin invariant} \label{part2}

Integrating the previous relations, we obtain an expression for the length of the caustic $\Gamma$ in terms of the Lazutkin invariant $L$ (recall that $\ell_0$ denotes the length of the boundary $\partial \Omega$):
{\small
\beqa{lengthGamma}
&&{\rm Length}(\Gamma) \;=\; \int_0^{2\pi} r(\f)\,d\f =\nonumber \\
&& \qquad =\; \int_0^{2\pi} \left[\rho(\f) - A(\f)\cdot L^{2/3} + B(\f)\cdot L^{4/3} + C(\f)\cdot L^{2} + D(\f)\cdot L^{8/3} + O\left(L^{10/3}\right)\right] d\f \nonumber\\
&& \qquad =\; \ell_0   - \left( \int_0^{2\pi} A(\f) \,d\f\right) \cdot L^{2/3} +  \left( \int_0^{2\pi} B(\f) \,d\f\right) \cdot L^{4/3} +  \left( \int_0^{2\pi} C(\f) \,d\f\right) \cdot L^{2} \nonumber\\
&&\qquad \quad +\;  \left( \int_0^{2\pi} D(\f) \,d\f\right) \cdot L^{8/3} + O\left(L^{10/3}\right)  \nonumber\\
&& \qquad =:\;
\ell_0   -  a \ L^{2/3} +  b \  L^{4/3} + c \  L^{2} 
+  d \  L^{8/3} + O\left(L^{10/3}\right)  
\eeqa}
where
\begin{eqnarray*}
a &=&     \frac{1}{2}  \left(\frac{3}{2}\right)^{2/3}  \left[\int_0^{2\pi} \rho^{1/3} \,d\f\right] \\
b &=&    \frac{1}{720} \left( \frac{3}{2}\right)^{1/3}  \left[
\int_0^{2\pi} \left(\frac{9}{\rho^{1/3}} +  8 \frac{\rho'^2}{\rho^{7/3}} \right)d\f \right]\\
c &=& 
\frac{1}{11200}  \left[ \int_0^{2\pi} \left(
\frac{9}{\rho} + 24 \left(\frac{\rho'^2 + \rho''^2 }{\rho^3}\right)
- 40 \frac{ \rho'^4}{\rho^5}\right) d\f \right] \\
d &=& 
\frac{1}{90} \left( \frac{3}{2}\right)^{2/3} 
\left[
\int_0^{2\pi} \left(
\frac{281}{44800}\cdot \frac{1}{\rho^{5/3}} +
\frac{1}{\rho^{11/3}} \left(
\frac{281}{8400}\rho'^2
+ \frac{167}{4200} \rho''^2
+
\frac{{\rho^{(3)}}^2}{42}
\right) 
+\frac{473 }{4725}\cdot
\frac{\rho''^3}{\rho^{14/3}} \right.\right. \\
&& \qquad \qquad \qquad \qquad  \left.\left.
- \; \frac{1}{\rho^{17/3}} \left(
\frac{11}{120} \rho'^4 + 
\frac{473}{945} \rho'^2 \rho''^2
\right)
+ \frac{781}{1458} \frac{\rho'^6}{\rho^{23/3}}
\right)d\f \right].
\end{eqnarray*}


\subsection{Rotation number of the caustic as a function of the  curvature of the boundary and the length of the caustic} \label{part3}

Let us denote  $I= - {\rm Length}(\Gamma)$. Oberve that $I \geq -\ell_0$ and it is equal to $-\ell_0$ when $L=0$ (on the boundary).  Now we would like to invert relation (\ref{lengthGamma}) to obtain an expansion of the Lazutkin invariant in terms of the length of the caustic (it plays the r\^ole of a {\it cohomology class}):
 {\small
\begin{eqnarray} \label{expralpha}
L(I) &=& \frac{(I+\ell_0)^{3/2}}{a^{3/2}}+ \frac{3\ b
  }{2\ a^{7/2}}  (I+\ell_0)^{5/2} + \frac{3
    \left(9\
   b^2+ 4\ a c\right)}{8\ a^{11/2}} (I+\ell_0)^{7/2} \nonumber \\
   && \;
   +\; \frac{ \left(24\ a^2 d+132\ a b c +143\
   b^3\right)}{16\ a^{15/2}} (I+\ell_0)^{9/2} + O\left((I+\ell_0)^{11/2}\right).
\end{eqnarray}
}

\vspace{10 pt}

This function corresponds to {\it Mather's $\alpha$ function} (at least for values of $I$ near $\ell_0$ for which there exists a caustic):
\begin{eqnarray*}
\a: [-\ell_0, -\ell_0+\delta) &\longrightarrow& \R\\
I &\longmapsto& L(I).\\
\end{eqnarray*}

This allows us to find  the rotation vector corresponding to a caustic $\G$ with Lazutkin invariant $L$ and length $-I$; recall, in fact, that
$\omega = \partial \alpha (I)$ (see subsection 1.3).

Therefore:

\begin{eqnarray} \label{omegaintermsofI}
\omega &=& \partial \alpha(I) = 
\frac{3 }{2
   a^{3/2}} \ (I+\ell_0)^{1/2} 
   + \frac{15 b}{4
   a^{7/2}} \ (I+\ell_0) ^{3/2} + \frac{21 \left( 9b^2 +4
   a c\right)}{16 a^{11/2}} \ (I+\ell_0)^{5/2} \nonumber \\
&& \qquad \qquad   +\; \frac{9 \left(24 a^2 d+132 a b
   c+143 b^3\right)}{32
   a^{15/2}} \ (I+\ell_0) ^{7/2}
   + O\left((I+\ell_0) ^{9/2}\right).
\end{eqnarray}

\subsection{Computing Mather's $\beta$-function on caustics} \label{part4}

Inverting the above expression, we obtain:
\begin{eqnarray*}
I(\omega) &=&  - \ell_0 +
\frac{4 a^3 }{9} \omega ^2 - \frac{80}{81}
   \left(a^4 b\right)  \omega ^4 +
   \left[\frac{16}{729} a^5 \left(87
   b^2-28 a c\right)+\frac{400 a^5
   b^2}{729}\right] \omega ^6 + \\
&& \;+ \;   \left[\frac{32 a^6 \left(-72 a^2 d+724
   a b c-909 b^3\right)}{6561}-\frac{160
   a^6 b \left(87 b^2-28 a
   c\right)}{6561}\right] \omega ^8 +O\left(\omega
   ^{10}\right).
\end{eqnarray*}

In conclusion, we obtain a representation of Mather's $\beta$-function at $\omega$:

\begin{eqnarray}
\beta(\omega) &=& I(\omega)\cdot \omega - \alpha\left( I(\omega) \right) = \label{betacaustic}\\
&=& 
-\ell_0  \omega +\frac{4 a^3 }{27} \omega
   ^3 -\frac{16}{81} 
   \left(a^4 b\right)\omega ^5 + \frac{64}{729}
   \left(4 a^5
   b^2- a^6 c\right) \omega ^7  \nonumber\\
   && \; - \; \frac{256 \left(3
   a^8 d-36 a^7 b c+56 a^6
   b^3\right)}{19683} \omega ^9  +O\left(\omega
   ^{11}\right). \nonumber
\end{eqnarray}

\subsection{Existence of caustics and end of the proof of Theorem \ref{maintheorem}} \label{part5}

In order to conclude the proof of Theorem \ref{maintheorem}, we need to address the following question: which billiards possess caustics?
We have already mentioned a negative result by John Mather \cite{Mather82} which says that 
caustics do not exist as soon as the curvature of the boundary vanishes at some point.\\
However, in our case -- {\it i.e.}, for strictly convex billiards -- the situation turns out to be completely different.

Let us recall an important 
result in the theory of billiards: Birkhoff billiards are {\it nearly-integrable}. In fact, in \cite{Lazutkin} V. Lazutkin 
introduced a very special change of coordinates that reduces 
the billiard map $f$  to a very simple form.

Let $L_\Omega:  {[0,\ell]\times [0,\pi] \to  \T \times [0,\dt]}$  with small 
$\dt>0$ be given by
\begin{equation}\label{Lazutkincoord}
L_\Omega(s,\phi)=\left(x=C^{-1}_\Omega \int_0^s k^{{2/3}}(s)ds,
\qquad
y=4C_\Omega^{-1}k^{{-1/3}}(s)\ \sin \frac{\phi}{2} \right),
\end{equation}
where $C_\Omega := \int_0^{\ell} k^{{2/3}} (s)ds$ is sometimes 
called the {\it Lazutkin perimeter} (observe that it is chosen so 
that period of $x$ is one).

In these new coordinates the billiard map becomes very simple (see \cite{Lazutkin}):

\begin{equation} \label{lazutkin-billiard-map}
f_L(x,y) = \Big( x+y +O(y^3),y + O(y^4) \Big)
\end{equation}

In particular, near the boundary $\{\phi=0\} = \{y=0\}$, the billiard map $f_L$ reduces to
a small perturbation of the integrable map $(x,y)\longmapsto (x+y,y)$, with a perturbation of size $O(y^3)$.

Using this result and an adapted version of KAM theorem, Lazutkin proved in \cite{Lazutkin}
that if $\partial \Omega$ is sufficiently smooth (smoothness is needed and determined by KAM theorem),
then there exists a positive measure set of caustics, which accumulates on
the boundary and on which the motion is smoothly conjugate to a rigid rotation
(see \cite{KP} for an improved version of Lazutkin's result). The corresponding rotation numbers form a positive measure Cantor set in the space of rotation numbers, which accumulates to zero (these rotation numbers are of {\it Diophantine type}).\\

This fact and the above discussion complete the proof of Theorem \ref{maintheorem}. In fact, on this positive-measure set of rotation numbers for which caustics exists, the above expression for $\beta$ holds and this family accumulates at $\omega=0$. In particular, $\beta$ is $C^{\infty}$ on a Cantor set, in the sense of Whitney (see also P\"oschel \cite{Poschel}).


We can recover from expression (\ref{betacaustic}) Taylor's coefficients of $\beta$-function: $\beta(\omega) = \sum_{k=0}^{\infty}  \beta_k \frac{\omega^k}{k!}$. 
First of all, $\beta_{2k}=0$ for all $k$'s (in fact, $\beta$ can be extended to an even function w.r.t. $\omega$). Then,
let us introduce the following {\it invariants} ($s$ denotes the arc-length and by $\dot{}$ we mean the derivative w.r.t  $s$):
\begin{eqnarray*}
\I_1 &:=&  \int_0^{2\pi} \rho \, d\f =  \int_{0}^{\ell_0}ds  = \ell_0 \\
\I_3 &:=& \int_0^{2\pi} \rho^{1/3} d\f   = \int_0^{\ell_0} k^{2/3} ds \\
\I_5 &:=& \int_0^{2\pi}  \left(\frac{9}{\rho^{1/3}} +  8 \frac{\rho'^2}{\rho^{7/3}} \right)d\f  
\ = \ \int_0^{\ell_0}  \left(\frac{ 9 + 8 \ \dot{\rho}^2 }{\rho^{4/3}}  \right)ds \\
&=&  \int_0^{\ell_0}  \left( 9\ k^{4/3} +  \frac{8 \ \dot{k}^2}{k^{8/3}}  \right)ds \\
\I_7 &:=&  \int_0^{2\pi} \left[
\frac{9}{\rho} + 24 \left(\frac{\rho'^2 + \rho''^2 }{\rho^3}\right)
- 40 \frac{ \rho'^4}{\rho^5}\right] d\f \\
&=& 
\int_0^{2\pi} \left[\frac{9}{\rho^2} + \frac{24}{\rho^2}
\left(\dot{\rho}^2 + \rho^2 \ddot{\rho}^2 + 2 \rho \dot{\rho}^2\ddot{\rho}
\right) - \frac{16 \dot{ \ \rho}^4}{\rho^2}
\right] ds \\
&=&  \int_0^{\ell_0} \left(
9\ k^2 + \frac{24\ \dot{k}^2}{k^2} + \frac{24\ \ddot{k}^2}{k^4}
- \frac{144\ \dot{k}^2 \ddot{k}}{k^5}
+ \frac{176\ \dot{k}^4}{k^6} 
\right) ds \\
\end{eqnarray*}

\begin{eqnarray*}
\I_9 &:=&
\int_0^{2\pi} \left[
\frac{281}{44800}\cdot \frac{1}{\rho^{5/3}} +
\frac{1}{\rho^{11/3}} \left(
\frac{281}{8400}\rho'^2
+ \frac{167}{4200} \rho''^2
+
\frac{{\rho^{(3)}}^2}{42}
\right) 
+\frac{473 }{4725}\cdot
\frac{\rho''^3}{\rho^{14/3}} \right. \\
&& \qquad  \left.
- \; \frac{1}{\rho^{17/3}} \left(
\frac{11}{120} \rho'^4 + 
\frac{473}{945} \rho'^2 \rho''^2
\right)
+ \frac{781}{1458} \frac{\rho'^6}{\rho^{23/3}}
\right) d\f\\
&=& \int_0^{\ell_0} 
\left[
\rho^{- 8/3} \left(
\frac{281}{44800 } +
\frac{281 }{8400} \dot{\rho}^2
-\frac{109 }{2100 } \dot{\rho}^4
+\frac{20291 }{127575 } \dot{\rho}^6\right)
+\rho^{-5/3} \left(
\frac{167 }{2100 }\dot{\rho}^2 \ddot{\rho}
- \; \frac{2411}{4725 }  \dot{\rho}^4 \ddot{\rho} \right)
 \right.\\
&& \qquad \left.
+\;
\rho^{-2/3} \left(
\frac{167 }{4200 }\ddot{\rho}^2
+\frac{122 }{675 } \dot{\rho}^2 \ddot{\rho}^2 
+  \frac{1 }{21} \dot{\rho}^3\ \dddot{\rho}
\right)
+ \rho^{1/3} \left( \frac{473  }{4725} \ddot{\rho}^3
+\frac{4}{21}
  \dot{\rho} \ddot{\rho}\ \dddot{\rho} \right)
+
\frac{1}{42} \rho^{4/3} \ {\dddot{\rho}}^2
\right] ds\\
&=&
\int_0^{\ell_0}
\left[
\frac{281}{44800} k^{8/3}
+\frac{281\ \dot{k}^2}{8400\ k^{4/3}} +
\frac{167\ \ddot{k}^2}{4200\ k^{10/3}}
- \frac{167 \ \dot{k}^2 \ \ddot{k}}{700\ k^{13/3}}
+ \frac{{\dddot{k}}^2}{42\ k^{16/3}}
+\frac{559 \ \dot{k}^4}{2100\ k^{16/3}} \right.\\
&& \left. \;
-\; \frac{473 \ \ddot{k}^3}{4725\ k^{19/3}}
-\frac{10 \ \dddot{k}\ \dot{k} \ \ddot{k}}{21\ k^{19/3}}
+
\frac{5 \ \dddot{k} \ \dot{k}^3}{7\ k^{22/3}}
+ \frac{13142 \ \dot{k}^2\ \ddot{k}^2}{4725\ k^{22/3}}
-\frac{10777 \ \dot{k}^4 \ \ddot{k}}{1575\ k^{25/3}}
+\frac{521897 \ \dot{k}^6}{127575\
k^{28/3}}
 \right] ds.
\end{eqnarray*}

\vspace{20 pt}

In particular, we have:

\begin{eqnarray*}
a &=&     \frac{1}{2}  \left(\frac{3}{2}\right)^{2/3}  \ \I_3  \\
b &=&    \frac{1}{720} \left( \frac{3}{2}\right)^{1/3}  \ \I_5 \nonumber\\
c &=& 
\frac{1}{11200}  \ \I_7 \\
d &=& 
\frac{1}{90} \left( \frac{3}{2}\right)^{2/3} 
\ \I_9 \nonumber
\end{eqnarray*}

and therefore:

\begin{eqnarray}\label{esprbeta}
\beta_1 &=& - \I_1  \nonumber\\
\beta_3 &=&  \frac{1}{4}\ \I_3^3\nonumber \\
\beta_5 &=& - \frac{1}{144} \ \I_3^4\ \I_5 \\
\beta_7 &=& \frac{1}{320} \ \I_3^5 \left( \frac{14}{81} \I_5^2 - \I_3 \I_7  \right) \ = \
\frac{\I_3^5 \left(14\ 
   \I_5^2-81\ \I_3
   \I_7\right)}{25920}\nonumber\\
 \beta_9 &=& - 7\ \I_3^6 \left(  \I_3^2 \ 
   \I_9- \frac{1}{5600}
   \I_3\ \I_5 \
   \I_7+\frac{7 }{583200} \I_5^3 \right).\nonumber
   \end{eqnarray}

Moreover, from (\ref{expralpha}), recalling the definition of the coefficients $a,b,c$ and $d$, one obtains:
$$
\a ( c) = \a_0\cdot {(c+\ell_0)^{3/2}} + \a_1\cdot (c+\ell_0)^{5/2} + \a_2\cdot (c+\ell_0)^{7/2} 
   + \a_3\cdot (c+\ell_0)^{9/2} + O\left((c+\ell_0)^{11/2}\right),
$$

where:

\begin{eqnarray*}
\alpha_0 &=&   \frac{4\sqrt{2}}{3}  \I_3^{-3/2} \\
\alpha_1 &=& \frac{\sqrt{2}}{135} \I_3^{-7/2} \I_5\\
\alpha_2 &=& \frac{1}{56700 \sqrt{2}}\left(
\frac{72\, {\I_3} {\I_7}+7 \,{\I_5}^2}{
   {\I_3}^{11/2}}   \right)\\
\alpha_3 &=& \frac{1}{826686000 \sqrt{2}} \left(
 \frac{261273600\, {\I_3}^2 {\I_9}+ 21384\,  {\I_3}
   {\I_5} {\I_7}+1001\, {\I_5}^3}{
   {\I_3}^{15/2}}\right).\\
\end{eqnarray*}

\qed

\vspace{10 pt}

As a byproduct, one could also compute the rotation vector as a function of the Lazutkin invariant $L$. In fact, from (\ref{lengthGamma}), (\ref{omegaintermsofI}) and the relation $I=-{\rm Length}(\Gamma)$, one obtains:
\begin{eqnarray*}
\omega(L) &=& 
 \left( \frac{3}{2}\right)^{1/3} \frac{2}{\I_3}\  {L}^{1/3} 
 +
 \frac{\I_5}{90\ \I_3^3}  
  L
 +
\left(\frac{3}{2} \right)^{2/3} \frac{
   \left(243  \  \I_3 \I_7+ 14 \ 
   \I_5^2\right)}{340200 \ 
   \I_3^3} \ L^{5/3}\\
   && \;+\;
   \left(\frac{3}{2} \right)^{1/3}
   \frac{
   \left(5443200\  \I_3^2
   \I_9 + 243 \ \I_3 \I_5
   \I_7+7\ 
   \I_5^3\right)}{30618000\
   \I_3^4} \ L^{7/3}+O\left(L^3\right)
\end{eqnarray*}

\noindent and its inverse:

\begin{eqnarray} \label{Lomega}
L(\omega)  &=&
\left[\frac{\I_3^3 }{12}\right] \ \omega^3
-    
\left[\frac{  \I_3^4
   \I_5}{4320}\right]\ \omega ^5
   +
\left[   \frac{\I_3^5  \left(14 \ \I_5^2-81 \ 
   \I_3
   \I_7 \right)}{21772800}
   \right]
   \ \omega ^7 \\
   && \;-\;
   \left[\frac{ \I_3^6 \left(4082400\ 
   \I_3^2 \I_9-729\  \I_3
   \I_5 \I_7+49\ 
   \I_5^3\right)}{26453952000
   } \right] \ \omega ^9+O\left(\omega ^{11}\right). \nonumber
\end{eqnarray}
Observe that this latter expression could be also obtained as
$L(\omega) = \omega \beta'(\omega) - \beta(\omega)$.\\


\subsection{Proof of Corollary \ref{maincor}} \label{part7}

Let us now prove Corollary \ref{maincor}. The proof easily follows from the expressions of $\b_1$ and $\b_3$, found in Theorem \ref{maintheorem}.
In fact, observe that:
$$
\beta_{3} + \pi^2 \beta_{1} \leq 0 \qquad
\Longleftrightarrow  \qquad 
\I_3^3 - 4\pi^2  \I_1 \leq 0.
$$
Now, using H\"older inequality (with $p=\frac{3}{2}$ and $q=3$):
\begin{eqnarray*}
\I_3 &=&  \int_0^{\ell_0} k^{2/3} ds  \leq \left( \int_0^{\ell_0} (k^{2/3})^{3/2} ds \right)^{2/3} 
\left( \int_0^{\ell_0} 1^{3} ds \right)^{1/3} \\
&=& (2\pi)^{2/3}  {\ell_0}^{1/3} = 
 (4\pi^2 \I_1)^{1/3}.
\end{eqnarray*}

Moreover, equality holds if and only if it holds in H\"older inequality. This means that $k$ must be constant (and strictly positive) and therefore, the curve must be a circle. \qed



\vspace{10 pt}



\vspace{30 pt}

\section{Some examples} \label{examples}

\subsection{Billiard in a disc}\label{exampledisc}
As we have already recalled in the Introduction, 
the billiard in a disc is one of the easiest examples of billiards. Let  $\cD$ be a disc of radius $R$.
It follows from elementary arguments that at each reflection of the ball the angle of incidence is the same as  the previous angle of reflection. Therefore, 
the angle of reflection remains constant along the orbit.
If we denote by $s$  the arc-length parameter ({\it i.e.}, $s\in \sfrac{\R}{ {\small 2\pi R \Z}}$) and by $\theta \in (0,\pi/2]$ the angle of reflection, then
the billiard map  has a very simple form:
$$
f(s,\theta) = (s + 2R\, \theta,\; \theta).
$$

 Let us now compute the previous invariants in this case. 

Let us start by observing that the $\beta$-function is given by:
$$
\beta(\omega) = - 2R \ \sin \left (\pi\omega \right ).
$$

Let us verify this. 
First of all, it is easy to check it for orbits of rotation number $\omega_n=1/n$. These orbits coincide with regular $n$-gons inscribed in $\cD$. It is easy to compute that each side of these polygons has length equal to  $ 2 R \sin \frac{\pi}{n} $ and therefore the total perimeter is $ 2n R \sin \frac{\pi}{n}$. It follows that
$$
\beta(1/n) = - \frac{1}{n} \left(   2n R \sin \frac{\pi}{n} \right) = -2 R \sin \frac{\pi}{n}. 
$$

More generally, the orbits of rotation number $\omega$ have (constant) angle of reflection $\theta=\pi \omega$  (it follow from the definition of rotation number and the fact that it must remain constant).
The segment joining two subsequent bounces have length $R \sin (\pi \omega)$ (see figure \ref{betadisc}), therefore it follows from the definition of $\beta$  (see (\ref{avaction})) that:
$$
\beta(\omega) =-  \lim_{N\rightarrow +\infty} \frac{1}{2N} \sum_{i=-N}^{N-1} R \sin (\pi \omega) = - R \sin (\pi \omega).
$$

\begin{figure} [h!]
\begin{center}
\includegraphics[scale=0.3]{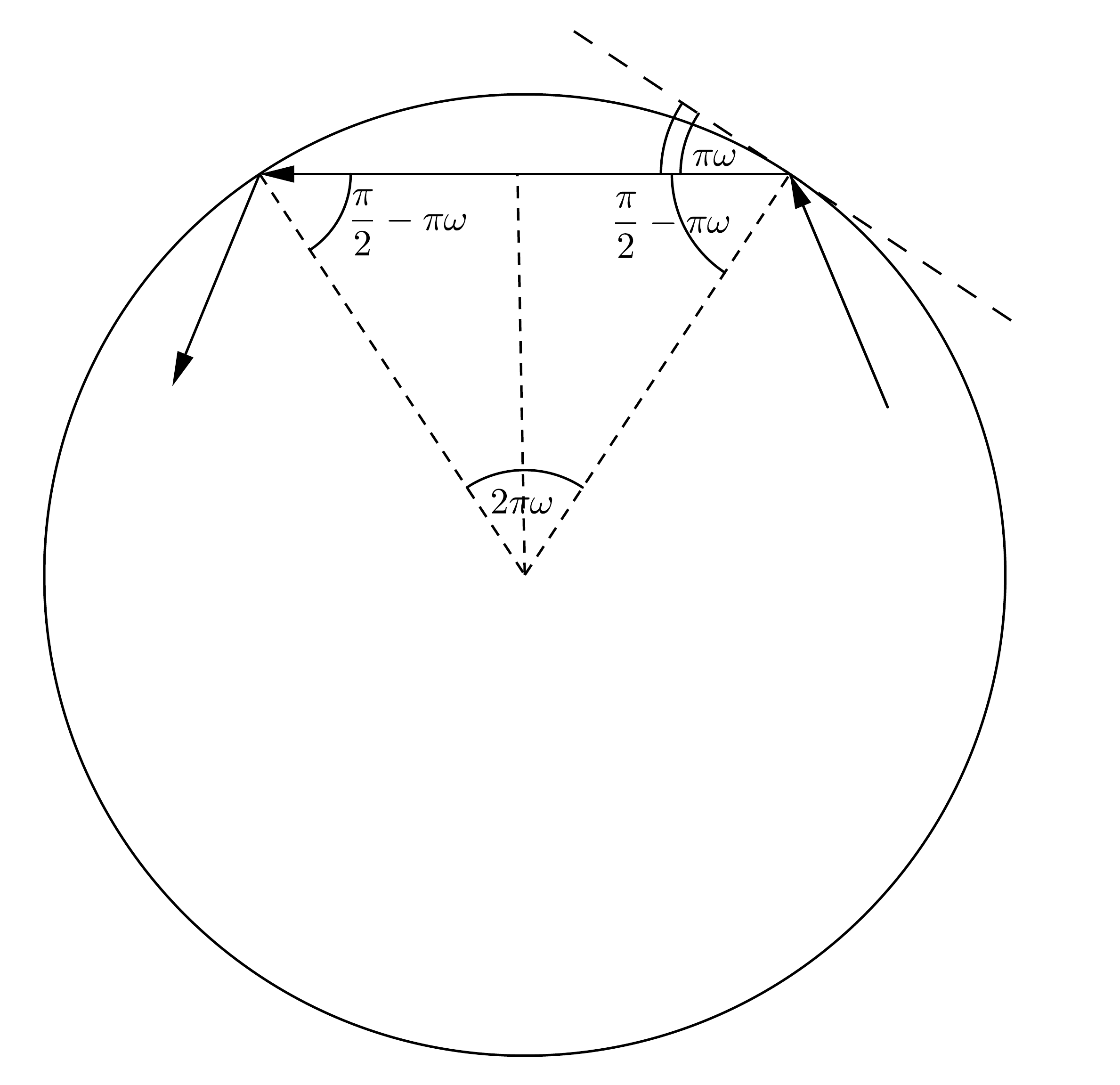}
\caption{}
\label{betadisc}
\end{center}
\end{figure}

Let us compute its Taylor expansion:

\begin{eqnarray}\label{betacerchio}
\beta(\omega) &=& - 2R \ \sin \left (\pi\omega \right ) \\
&= &
- 2 \pi  R \ \omega +\frac{1}{3!} \left(2R \pi ^3
   \right) \omega^3+\frac{1}{5!} \left (- 2 R \pi^5
   \right) \omega^5  + 
   \frac{1}{7!} \left(2R \pi^7 \right)  \omega^7 + \frac{1}{9!} \left( - 2R \pi^9\right)
   \omega^9+O\left(\omega^{11}\right
   ). \nonumber
\end{eqnarray}

In particular, $k(s)\equiv \frac{1}{R}$  and $\ell_0=2\pi R$. The above invariants are therefore:
\begin{eqnarray}
\I_1 &=& 2\pi R \nonumber\\
\I_3 &=& 2\pi R^{1/3} \nonumber\\
\I_5 &=& 18 \pi R^{-1/3}  \label{invcerchio}\\
\I_7 &=& 18 \pi R^{-1} \nonumber\\
\I_9 &=& \frac{281}{22400} \pi R^{-5/3}.\nonumber
\end{eqnarray}

Substituting these values in (\ref{esprbeta}), one can easily check that they match with (\ref{betacerchio}).

Moreover, one also obtains (by geometric reasoning) that:
\begin{eqnarray*}
L(\omega) &=& 
-2\pi R \omega \cos (\pi \omega) + 2R \sin (\pi \omega) \\
&=& 
\frac{2}{3} \pi ^3 R \ \omega^3-\frac{1}{15} \left(\pi ^5 R\right)  \omega ^5+\frac{1}{420} \pi ^7 R \ \omega^7-
\frac{\left(\pi^9 R\right) }{22680} \omega^9 +O\left(\omega ^{11}\right).
\end{eqnarray*}
One can check that this expression matches with  (\ref{Lomega}) once the invariants (\ref{invcerchio}) are substituted in it.

Finally, observe that:
\begin{eqnarray*}
\omega(L) &=& 
\left( \frac{3}{2}\right)^{1/3}
\frac{{L}^{1/3}}{\pi  {R^{1/3}}}+\frac{L}{20\ \pi 
   R}+
   \left( \frac{3}{2}\right)^{2/3}
   \frac{41 \ L^{5/3}}{8400\
 \pi  \ R^{5/3}}+
 \left( \frac{3}{2}\right)^{1/3}
 \frac{97\ 
   L^{7/3}}{100800\  \pi 
   \ R^{7/3}}+O\left(L^3\right).\\
   \end{eqnarray*}

%
%
%
%
%

\vspace{10 pt}

\subsection{Billiard in an ellipse}\label{exampleellipse}

Let us consider now the billiar inside an ellipse 
$$
\cE = \left\{(x,y): \;\frac{x^2}{a^2}+ \frac{y^2}{b^2}=1\right\}
$$
with $0<b\leq a$. Up to rescaling, we can assume that $a=1$ (see also Remark \ref{ingeneral}) and therefore the eccentricity of the ellipse is given by
$0\leq h=\sqrt{1-b^2}<1$ and the two foci by $F_{\pm}= (\pm h, 0)$.

Optical properties of conics (an alternative way to consider the billiard ball motion inside a conic) were already well known to ancient Greeks. We refer to \cite{Tabach} for a more detailed discussion (see also \cite{Siburg}). In particular, billiard trajectories can be classified in the following way:
\begin{itemize}
\item[a)] trajectories that always intersect the open segment between the two foci,
\item[b)] trajectories that never intersect the closed segment between the two foci, and
\item[c)] trajectories that alternatively pass through one of the two foci.
\end{itemize}

In particular, each trajectory in a) is tangent to a confocal hyperbola, each trajectory in b) is tangent to a confocal ellipse, while trajectories of kind c) tend asymptotically to the major semiaxis.
Confocal ellipses are therefore examples of caustics (also hyperbolae can be considered a sort of  generalized caustics) which foliate everything but the closed segment between the two foci (see figure \ref{ellipse-billiard} in subsection 1.2). Hence, this is an example of an integrable billiard, as we have already recalled in the Introduction.\\

Let us now try to describe the dynamics and provide some expression for its $\beta$-function.  Differently from the circular case, here the situation is much more complicated due to the appearance of elliptic integrals, which make the dynamics much less explicit. A description of the dynamics is carried out, for example, in \cite{Taba}.

Let us introduce the following elliptic coordinates
$$
\cE_{\m}: \quad \left\{
\begin{array}{l}
x= h \cosh \m \, \cos \f\\
y= h \sinh \m \, \sin \f\\
\end{array}\right. \qquad \f\in [0,2\pi),\; 0<\m\leq \m_0,
$$
where $\m_0$ is such that $\cosh \mu_0 = \frac{1}{h}$. Observe that $\cE_{\m_0}$ corresponds to our boundary ellipses, while $\cE_{\m}$ are the 
confocal ones.

Let us denote $I(\m)=\cosh^2 \m$. In particular, the lengths of these caustics are (let us denote by $a_\mu$ the major semi-axis of $\cE_\mu$):
$$
|\cE_{\m}| =  4 a_\mu E\left(\frac{1}{\sqrt{I(\mu)}}\right) =
4 h \sqrt{I(\mu)}\, E\left(\frac{1}{\sqrt{I(\mu)}}\right),
$$
where $E(k)=\int_0^{\pi/2} \sqrt{1-k^2\sin^2\theta}d\theta$ is a complete elliptic integral of second type.\\

It follows from \cite[formula 1.7]{Taba} that\footnote{The different factor in front of it, follows from a slightly different definition of rotation number.}:
$$
\omega(I(\mu)) = {\rm Rot}(I(\mu)) = \frac{1}{4F\left(\frac{1}{\sqrt{I}}\right)}
F\left(\arcsin \left( \sqrt{I} \frac{2\tanh \m_0 \sqrt{\cosh^2\m_0-I}}{\cosh^2\m_0-I+I\tanh^2\m_0}
\right), \frac{1}{\sqrt{I}}
\right)
$$
where $F(z,k)=\int_0^z \frac{d\theta}{\sqrt{1-k^2\sin^2\theta}}$ denotes an elliptic integral of first type. In the following, we shall denote the complete elliptic integral of
first type by $K(k)=F(\frac{\pi}{2},k)$.

Using these results we can compute Mather's $\beta$-function in this case. Here are the needed steps:
\begin{itemize}
\item[I - ] let $c(I)= - |\cE_{\mu}|$ denote the cohomology class, seen as a function of $I$;
 \item[II - ] one could invert the function $\omega(I)$, which is a function of $\sqrt{I-\frac{1}{h^2}}$, and obtain a function
 $I(\omega)$ in a neighbourhood of $\omega=0$;
 \item[III - ] then, one obtains an expression $c=c(\omega)$; recalling that $c(\omega)=\beta'(\omega)$ and integrating, one finds an expression for $\beta$.
  \end{itemize}
  
Carrying out these computations, we get:
\begin{eqnarray} \label{betaellipse}
  \beta(\omega) &=&
   - 4E(h) \, \omega 
   + \left[\frac{8}{3}(1-h^2) K(h)^3\right] \omega^3 
   +  \left[\frac{8}{15} (1-h^2) K(h)^4 \left[  15 E(h) - 8(2-h^2) K(h) \right]\right] \,\omega^5 \nonumber \\
&+& \left[   \frac{16}{315} \left(1-h^2\right) K(h)^5 
\left[ 630 E(h)^2 - 630 \left(2- h^2\right) K(h) E(h)+\left(136 h^4-631
   h^2+631\right) K(h)^2 \right]\right]\omega^7\nonumber \\
   &-&
    \left[
    \frac{8 \left(h^2-1\right) K(h)^6 \left[75600
   \left(h^2-2\right) K(h) E(h)^2+4
   \left(h^2-2\right) \left(992 h^4-5741 h^2+5741\right)
   K(h)^3\right.}{2835} \right.\nonumber\\
   && \qquad \left. \frac{\left. + 1323 \left(24 h^4-109 h^2+109\right)
   K(h)^2 E(h)+52920
   E(h)^3\right]}{2835}
   \right]
   \omega^9
   + O(\omega^{11}).
\end{eqnarray}
    
    It is easy to check that in the limit as $h\rightarrow 0$, we recover the $\beta$-function for the circular billiard of radius $R=1$. Observe in fact that:
    $$
    \lim_{h\rightarrow 0^+} E(h) = 
        \lim_{h\rightarrow 0^+} K(h) = \frac{\pi}{2}.$$
        
\vspace{10 pt}
      
  We can also verify this expression, computing the invariants $\I_k$'s directly from Theorem \ref{maintheorem}. For the sake of this presentation, we shall compute only  $\I_1, \I_3$ and $\I_5$ and verify the corresponding coefficients $\b_1, \b_2, \b_3$. The others could be computed similarly, but, for  simplicity, we omit those lenghty -- yet, similar -- computations. 
  
  Let us consider the parametrizion of $\cE_{\mu_0}$ by polar coordinates (as above).
  Recall that the arc-length is given by 
  $ds= \sqrt{1-h^2\cos^2 \f}\, d\f$, while  the curvature in polar coordinate is:
  $$
  k(\f) = \frac{\sqrt{1-h^2}}{(1-h^2\cos^2 \f)^{3/2}}.
  $$
  It is easy to check that:
  \begin{itemize}
  \item[i)] 
  \begin{eqnarray*}
  \I_1 &=& \int_0^{\ell_0} ds = \int_0^{2\pi} \sqrt{1-h^2\cos^2 \f} \,d\f = 4  \int_0^{\frac{\pi}{2}} \sqrt{1-h^2\sin^2 \f} \,d\f = 4E(h).
  \end{eqnarray*}
  Therefore, $\b_1= - 4E(h)$, which matches with the expression in (\ref{betaellipse}).
    \item[ii)]
      \begin{eqnarray*}
  \I_3 &=& \int_0^{\ell_0} k^{2/3} ds = 
  \int_0^{2\pi}  \frac{(1-h^2)^{1/3}}{1-h^2\cos^2\f}  \sqrt{1-h^2\cos^2 \f} \, d\f =\\
  &=&
4 (1-h^2)^{1/3}  \int_0^{\frac{\pi}{2}}  \frac{1}{\sqrt{1-h^2\cos^2\f}}  \, d\f = 4 (1-h^2)^{1/3} K(h).\\  
  \end{eqnarray*}
  In particular, using Theorem \ref{maintheorem} we also obtain
  \begin{eqnarray*}
  \frac{\b_3}{3!} &=& \frac{1}{3!} \left( \frac{1}{4}\I_3^3\right) = \frac{8}{3} (1-h^2) K(h)^3.
  \end{eqnarray*}
 \item[iii)]  First of all, let us observe that 
 \begin{eqnarray*}
 \dot{k} &=& \frac{dk}{ds} = \frac{dk}{d\f} \cdot \frac{d\f}{ds} =\\
 &=& \frac{dk}{d\f} \cdot  \frac{1}{\sqrt{1-h^2\cos^2\f}} =\\
 &=& - 3 h^2 \sqrt{1-h^2} \frac{\cos \f \sin\f}{(1-h^2 \cos^2\f)^3}.
 \end{eqnarray*}
Hence:
     \begin{eqnarray*}
  \I_5 &=&\int_0^{\ell_0}  \left( 9\ k^{4/3} +  \frac{8 \ \dot{k}^2}{k^{8/3}}  \right)ds =\\
  &=& \int_0^{2\pi}  \left( 
  9 \frac{(1-h^2)^{2/3}}{(1-h^2\cos^2\f)^{3/2}}
   + 
     72 \, h^4 (1-h^2)^{-1/3} \frac{\cos^2\f \sin^2\f}{(1-h^2\cos^2\f)^{3/2}}
       \right) d\f  =\\
    &=&  4 \int_0^{\frac{\pi}{2}}  \left( 9 \frac{(1-h^2)^{2/3}}{(1-h^2\cos^2\f)^{3/2}}
    + 
     72 \, h^4 (1-h^2)^{-1/3} \frac{\cos^2\f \sin^2\f}{(1-h^2\cos^2\f)^{3/2}}
       \right)d\f  =\\   
       &=& 36 (1-h^2)^{-1/3} E(h) + 288 (1-h^2)^{-1/3} \left(  - 2E(h) + (2-h^2) K(h)\right) =\\
       &=& 36 (1-h^2)^{-1/3}  \left[  - 15E(h) + 8 (2-h^2) K(h)\right].\\       
  \end{eqnarray*}
Substituting in the expression of $\b_5$ (see Theorem \ref{maintheorem}):
\begin{eqnarray*}
\frac{\b_5}{5!} &=& -\frac{1}{144 \cdot 5!} \I_3^4 \I_5 =\\
&=& -\frac{4^4}{144 \cdot 5!}  \, (1-h^2)^{4/3} K(h)^4    \left[ 36 (1-h^2)^{-1/3}  \left(  - 15E(h) + 8 (2-h^2) K(h)\right)\right] =\\
&=& \frac{8}{15} (1-h^2) K(h)^4  \left[  15E(h) - 8 (2-h^2) K(h)\right].\\
\end{eqnarray*}
  \end{itemize}
        
 In the same way one could compute $\I_7$ and $\I_9$.\\

 \begin{Rem}\label{ingeneral}
 Similar formulae hold in the general case, {\it i.e.}, without assuming that the major semiaxis $a=1$. Let us consider an ellipse $\cE$ with semiaxis
 $0<b\leq a$ and eccentricity $h=\sqrt{1-\left(\frac{b}{a}\right)^2}$. It follows easily from the definition of $\beta$-function, that rescaling the ellipse, this function will rescale by the same amount. Therefore, one could consider the rescaled ellipse $\frac{1}{a}\cE$ -- which has major semiaxis equal to 1 and the same eccentricity $h$ as $\cE$ -- and use the above formulae for computing the corresponding $\beta$-function. The $\beta$-function associated to the original ellipse $\cE$ will be given by
$
\beta_{\cE} = a\beta_{\frac{1}{a}\cE}.
$
 \end{Rem}

 To conclude this section, we would like to address the following question: is it true that the $\beta$-function determines univocally an ellipse amongst other ellipses? In other words: is it possible that two different ellipses have the same $\beta$-function? We shall show that the first question (resp. the second question) has an affirmative answer (resp. negative answer).
 
 \begin{Prop}\label{uniqueness}
 If $\cE_1$ and $\cE_2$ are two ellipses such that $\beta_{\cE_1} \equiv \beta_{\cE_2}$, then $\cE_1$ and $\cE_2$ are the same ellipse. More generally: if the Taylor coefficients 
$ \beta_{\cE_1,1} =  \beta_{\cE_2,1}$ and $ \beta_{\cE_1,3} =  \beta_{\cE_2,3}$, then the same conclusion remains  true.\\
 \end{Prop}

\begin{Proof}
We prove the second statement, which clearly implies the first one. Let us denote by $0<b_i\leq a_i$ the semi-axis of $\cE_i$, and by $h_i=\sqrt{1-\left(\frac{b_i}{a_i}\right)^2}$ their eccentricities.
If $ \beta_{\cE_1,1} =  \beta_{\cE_2,1}$ and $ \beta_{\cE_1,3} =  \beta_{\cE_2,3}$, then using the above expressions and Remark \ref{ingeneral}, we can conclude that:
\begin{equation}\label{sistema}
\left\{
\begin{array}{l}
a_1 \,E(h_1) = a_2\, E(h_2) \\
a_1 \, (1-h_1^2)\, K(h_1)^3 = a_2 \, (1-h_2^2)\, K(h_2)^3.\\
\end{array}
\right.
\end{equation}
In particular, since $a_i\neq 0$ and $E(h_i)\neq 0$, it follows that:

$$
\frac{(1-h_1^2)\, K(h_1)^3 }{E(h_1)}= \frac{ (1-h_2^2)\, K(h_2)^3}{E(h_2)}.
$$

One can check that the function $f(x)=\frac{(1-x^2)\, K(x)^3 }{E(x)}$ is strictly decreasing\footnote{
One could also compute its first derivative explicitly and show that it is stricly negative in $(0,1)$ and tends to $-\infty$ in the limit as $x\rightarrow 1^-$:
$$
f'(x) = \frac{K(x)^2}{x\, E(x)^2} \Big( 3E(x)^2 + 2(2-x^2)E(x) K(x) + (1-x^2) K(x)^2\Big).
$$
} in $[0,1]$, with  $f(0)= \frac{\pi^2}{4} $  (which corresponds to the circular case) and $f(1)= 0$ (degeneration of the ellipse into a parabola). 
Therefore, if $f(h_1)=f(h_2)$, then $h_1=h_2$, {\it i.e.},  the two ellipses have the same eccentricity. Substituting this piece of information in the first equation of (\ref{sistema}), one also obtains that 
$a_1=a_2$ and consequently $b_1=b_2$. This concludes the proof.
\end{Proof}


\vspace{1.truecm}

\end{document}